\documentclass{article}

\usepackage{graphicx}
\usepackage[german,english]{babel}
\usepackage{amsmath}
\usepackage{amssymb}
\usepackage{amsfonts}
\usepackage{bbm}
\usepackage{url}
\usepackage{booktabs}
\usepackage{tabularx}
\usepackage{epsfig,subfigure}
\usepackage{longtable}
\usepackage{lscape}
\usepackage{psfrag}

\newtheorem{defn}{Definition}

\newtheorem{thm}[defn]{Theorem}
\newtheorem{cor}[defn]{Corollary}
\newtheorem{conj}[defn]{Conjecture}

\makeatletter

\graphicspath{{.}}

\title{Triangulated Manifolds with Few Vertices: Centrally Symmetric Spheres and Products of Spheres}

\author{\Large Frank H.~Lutz}
\date{}

\begin{document}

\selectlanguage{english}

\maketitle

\bigskip
\bigskip
\bigskip

\noindent
Let $M$ be a simplicial 
manifold with $n$ vertices.
We call $M$ \emph{centrally symmetric} if it is invariant under 
an involution $I$ of its vertex set which fixes no face of $M$.
Obviously, the number of vertices of a centrally symmetric
(triangulated) manifold is even, $n=2k$, and, 
without loss of generality, we may assume that the involution 
is presented by the permutation $I=(1\,\,\,k+1)(2\,\,\,k+2)\cdots (k\,\,\,2k)$.\linebreak
The boundary complex $\partial C_k^{\Delta}$ of the $k$-dimensional
crosspolytope $C_k^{\Delta}$ is clearly centrally symmetric with
respect to the standard antipodal action, 
and a subset $F\subseteq\{ 1,2,\ldots ,2k\}$ 
is a face of $\partial C_k^{\Delta}$ if and only if 
it does not contain any \emph{minimal non-face} $\{ i,k+i\}$ for $1\leq i\leq k$.
Hence, every centrally symmetric manifold with 
$2k$ vertices appears as a subcomplex of the boundary complex 
of the $k$-dimensional crosspolytope. 

Free ${\mathbb Z}_2$-actions on spheres are at the heart 
of the Borsuk-Ulam theorem, which has an abundance 
of applications in topology, combinatorics, functional analysis,
and other areas of mathematics (see the surveys of Steinlein
\cite{Steinlein1985}, \cite{Steinlein1993}, 
and the recent book of Matou\v{s}ek \cite{Matousek2003}).
Centrally symmetric spheres therefore constitute an important 
class of triangulated spheres for which we have a strong interest 
in understanding their combinatorial properties, like the range 
of possible $f$-vectors, or even more basic, what kind of examples 
are there at all? 

Centrally symmetric products of spheres are the next more 
general class of centrally symmetric manifolds.
They 
show that certain
lower bounds on the numbers of vertices of 
centrally symmetric manifolds are tight.

The aim of this 
paper
is to give a survey of the known results
concerning centrally symmetric polytopes, spheres, and manifolds.
We further enumerate \emph{nearly neighborly} centrally symmetric spheres and
centrally symmetric products of spheres with dihedral or cyclic symmetry on few vertices,
and we present an infinite series of vertex-transitive nearly neighborly 
centrally symmetric $3$-spheres.

\subsection*{1\,\, General Properties of Centrally Symmetric Spheres}

One way to obtain centrally symmetric spheres is as boundary
complexes of centrally symmetric simplicial polytopes. 
A $d$-dimensional polytope $P\subset {\mathbb R^{d}}$ 
is \emph{centrally symmetric}
if we can translate $P$ such that $P=-P$. If $d>0$, then, by convexity, the involution
$I:x\mapsto -x$ of $\mathbb R^{d}$ does not fix any non-trivial face of $P$, 
and $P$ has an even number of vertices, $n=2k$.
Regular $2k$-gons, the icosahedron, and crosspolytopes $C_k^{\Delta}$
are immediate examples of centrally symmetric simplicial polytopes.
The dodecahedron and $d$-dimensional cubes are centrally symmetric, 
but not simplicial. 

Not every centrally symmetric sphere needs to be polytopal, 
and even if so, resulting realizations need not be centrally
symmetric. Centrally symmetric simplicial $(d-1)$-spheres have at least $2d$
vertices, with the boundary complex $\partial C_{d}^{\Delta}$
of the $d$-dimensional crosspolytope $C_{d}^{\Delta}$
as the unique centrally symmetric $(d-1)$-sphere with exactly
$2d$ vertices.

We recall that for the class of \emph{all} simplicial spheres,
the upper bound theorem of McMullen \cite{McMullen1970} for polytopal spheres
and of Stanley~\cite{Stanley1975} for simplicial spheres
(see Novik \cite{Novik1998} for generalizations to 
odd-dimensional and certain even-dimensional simplicial manifolds)
as well as the lower bound theorem of Barnette
(\cite[p.~354]{Barnette1973a}, \cite{Barnette1973b}) and Kalai \cite{Kalai1987}
give restrictions on the numbers $f_i$ of $i$-dimensional faces of a
simplicial sphere for $0\leq i\leq d-1$:
A simplicial $(d-1)$-sphere with $n$ vertices has at most as many $i$-faces
as the boundary sphere of the corresponding cyclic $d$-polytope $C_{d}(n)$
and at least as many $i$-faces as the boundary sphere of a stacked $d$-poly\-tope
on $n$ vertices. In contrast, much less is known on $f$-vectors $f=(f_0,\dots,f_{d-1})$
of centrally symmetric $d$-polytopes respectively $(d-1)$-spheres.

Stanley \cite{Stanley1987} proved lower bounds (conjectured by 
B\'ar\'any and Lov\'asz \cite{BaranyLovasz1982} and by Bj\"orner)
on the numbers of faces of $d$-dimensional centrally symmetric polytopes 
with $n=2k\geq 2d$ vertices (see Novik \cite{Novik1999} for an alternative 
and more geometric proof):
$$f_i\geq 2^{i+1}\binom{d}{i+1}+2(k-d)\binom{d}{i},\quad 0\leq i\leq d-2,$$
$$f_{d-1}\geq 2^d+2(k-d)(d-1).$$
These bounds are sharp for \emph{stacked centrally symmetric $d$-polytopes},
which are obtained from the $d$-dimensional crosspolytope
by stellarly subdividing $n-k$ successive pairs of 
antipodal facets.

A simplicial $(d-1)$-sphere $S$ is \emph{$l$-neighborly} if every set of $l$
(or less) vertices forms a face of $S$. The $d$-simplex $\Delta_{d}$ (respectively,
its boundary $\partial\Delta_{d}$)
with $d+1$ vertices is $(d+1)$-neighborly,
and for $n\geq d+2$, the cyclic polytope $C_{d}(n)$
is $\lfloor \frac{d}{2}\rfloor$-neighborly, 
but not $(\lfloor \frac{d}{2}\rfloor+1)$-neighborly.
Simplicial spheres (respectively, simplicial polytopes) are called \emph{neighborly}
if they are $\lfloor \frac{d}{2}\rfloor$-neighborly.

Analogously, a centrally symmetric $(d-1)$-sphere $S$ with $n=2k$ vertices 
is \emph{centrally $l$-neighborly}
if every set of $l$ vertices, which does not contain
a minimal non-face $\{ i,k+i\}$ for $1\leq i\leq k$, 
is a face of $S$, i.e., if $S$ has the $(l-1)$-skeleton 
of the crosspolytope $C_k^{\Delta}$.
The $d$-dimensional crosspolytope
$C_{d}^{\Delta}$ with $2d$ vertices
is centrally $d$-neighborly. 
A centrally symmetric $(d-1)$-sphere with $n=2k$ vertices is \emph{nearly neighborly}
if it is centrally $\lfloor\frac{d}{2}\rfloor$-neigh\-borly, i.e.,
if $f_i=2^{i+1}\binom{k}{i+1}$ for $i\leq\frac{d}{2}-1$,
with $f_i$ being determined by the Dehn-Sommerville equations
for $i >\frac{d}{2}-1$.

Along the lines of the proof of the upper bound theorem for simplicial spheres,  
Adin \cite{Adin1991} and Stanley (cf.~\cite{Jockusch1995})
showed independently that a centrally symmetric
simplicial $(d-1)$-sphere with $2k$ vertices has at most as many $i$-faces as a 
nearly neighborly centrally symmetric $(d-1)$-sphere with $2k$
vertices would have, if such exists. Novik \cite{Novik2003} extended this result to all odd-dimensional
centrally symmetric manifolds; see also \cite{Novik2003pre}.

The boundaries of regular polygons with $2k\geq 4$ vertices and suspensions 
thereof with $2k+2$ vertices provide examples of
centrally symmetric $1$- and $2$-spheres for all possible
numbers of vertices. Since centrally $1$-neighborliness is a trivial property,
every centrally symmetric $2$-sphere is nearly neighborly, 
and, moreover, is realizable as the boundary complex of a 
centrally symmetric $3$-polytope; see Mani \cite{Mani1971}.

Gr\"unbaum observed \cite[p.~116]{Gruenbaum1967} that the centrally symmetric $4$-polytope 
$G_{2\cdot 4+2}^{4}:=\mbox{conv}\{\pm e_1,\dots,\pm e_{4},\pm {\mathbbm 1}\}\subset {\mathbb R^{4}}$
on $2\cdot 4+2$ vertices is simplicial and nearly neighborly, 
but that there are \emph{no} nearly neighborly 
centrally symmetric $4$-polytopes with $n\geq 12=2\cdot 4+4$ vertices.
In fact, McMullen and Shephard~\cite{McMullenShephard1968} proved that 
centrally symmetric $d$-poly\-topes with $n\geq 2d+4$ vertices
are at most centrally $\lfloor\frac{d+1}{3}\rfloor$-neigh\-borly.
Hence, there are no nearly neighborly centrally symmetric $d$-poly\-topes 
with $n\geq 2d+4$ vertices \emph{for all} $d\geq 4$.
According to Pfeifle \cite[Ch.~10]{Pfeifle2003} also nearly neighborly centrally symmetric
$d$-dim\-en\-sional fans on $2d+4$ rays do not exist for all even $d\geq 4$ and all odd $d\geq 11$.
Schneider \cite{Schneider1975} gave an asymptotic lower bound for the maximal possible $l=l(d,s)$ 
for which there are centrally $l$-neighborly $d$-polytopes with $2(d+s)$ vertices.
However, Burton~\cite{Burton1991} showed that, for fixed dimension
$d\geq 4$, centrally symmetric $d$-polytopes with sufficiently many 
vertices \emph{cannot} be centrally $2$-neighborly.

In contrast to the situation for centrally symmetric polytopes,
Gr\"unbaum constructed nearly neighborly centrally symmetric $3$-spheres 
with $12$ and $14$ vertices; see \cite{Gruenbaum1969}, \cite{Gruenbaum1970c}, 
and \cite{Gruenbaum1970a}.

\bigskip

\noindent
\textbf{Centrally Symmetric Upper Bound Conjecture}
(Gr\"un\-baum~\cite{Gruenbaum1970a})\label{conj:gruenbaum}\linebreak
\emph{There are nearly neighborly centrally symmetric $(d-1)$-spheres 
with $n$ vertices for all\, $d\geq 2$ and even\, $n=2k\geq 2d$.}

\bigskip

Since being centrally $\lfloor\frac{d}{2}\rfloor$-neighborly is preserved 
under suspension and since $\lfloor\frac{d}{2}\rfloor=\lfloor\frac{d+1}{2}\rfloor$ for all even $d$,
it suffices to construct \emph{odd-dimensional} nearly neighborly
centrally symmetric $(d-1)$-spheres for all even numbers $n\geq 2d$ of vertices in order to verify
Gr\"unbaum's centrally symmetric upper bound conjecture.

Gr\"unbaum's conjecture is trivial for $1$- and $2$-spheres, but also
holds for $3$- and $4$-spheres.

\begin{thm} {\rm (Jockusch \cite{Jockusch1995})}
There is an infinite family\, $J^{\,3}_{2k}$, $k\geq 4$, of \linebreak
nearly neighborly centrally
symmetric $3$-spheres with\, $2k$ vertices. Moreover,
the suspensions $S^0*J^{\,3}_{2k}$ form a family of 
nearly neighborly centrally symmetric $4$-spheres with\, 
$2k+2$ vertices for $k\geq 4$.
\end{thm}

Jockusch constructs the series $J^{\,3}_{2k}$ by induction. 
He starts with the boundary complex $J^{\,3}_8=\partial C_{4}^{\Delta}$ 
of the $4$-dimensional crosspolytope with $8$ vertices. 
For the induction step he chooses a $3$-ball $B_{2k}$ 
with image $B_{2k}^I$ under the central symmetry $I$
such that their intersection\, $B_{2k}\cap\, B_{2k}^I$\,
does not contain any facet of\, $J^{\,3}_{2k}$.
He then removes the balls\, $B_{2k}$\, and\, $B_{2k}^I$\, 
from $J^{\,3}_{2k}$ and sews in two new balls\,
$(2k+1)*\partial B_{2k}$\, and\, $(2k+2)*\partial B_{2k}^I$\, 
to obtain the \mbox{$3$-sphere}\, $J^{\,3}_{2k+2}$.
The way Jockusch chooses the balls $B_{2k}$ (the balls $B_{2k}$ and $B_{2k}^I$ contain
all the vertices of $J^{\,3}_{2k}$, but have no interiour edges, respectively), 
he ensures that\, $J^{\,3}_{2k+2}$\,
remains centrally symmetric and nearly neighborly in every step. 

\begin{thm} {\rm (McMullen and Shephard \cite{McMullenShephard1968})}
For even $d$, let the polytope \linebreak $H_{\,2d+2}^{\,d}:={\rm conv}(\Delta_d\cup -\Delta_d)$ 
be the joint convex hull of a regular $d$-simplex $\Delta_d$ (with center $0$)
and its image $-\Delta_d$ under the map $I:x\mapsto -x$.
Then $H_{\,2d+2}^{\,d}$ is nearly neighborly and has the group $S_{d+1}\times{\mathbb Z}_2$
as its vertex-transitive geometric automorphism group. 
\end{thm}

Gr\"unbaum \cite[p.~116]{Gruenbaum1967} has shown that there is only one combinatorial type of a nearly neighborly centrally symmetric
$4$-polytope with $10$ vertices, i.e., $G_{2\cdot 4+2}^{4}$ and
$H_{\,2\cdot 4+2}^{\,4}$ are combinatorially isomorphic
(in fact, for all even $d$ $G_{2\cdot d+2}^{d}:=\mbox{conv}\{\pm e_1,\dots,\pm e_{d},\pm {\mathbbm 1}\}$ 
is combinatorially isomorphic to $H_{\,2d+2}^{\,4}$).

In odd dimensions $d+1$ the polytope $H_{\,2(d+1)+2}^{\,d+1}$ is not simplicial.
However, ${\rm conv}((\Delta_d\cup -\Delta_d)\cup \{\pm e_{d+1}\})\subset {\mathbb R^{d+1}}$
is a nearly neighborly centrally symmetric $(d+1)$-dimensional polytope on $2d+4$ vertices
with boundary \linebreak
$\partial\, {\rm conv}((\Delta_d\cup -\Delta_d)\cup \{\pm e_{d+1}\})=S^0*H_{\,2d+2}^{\,d}$.

If $d$ is even, then, on the combinatorial level, 
the sphere $\partial H_{\,2d+2}^{\,d}$
can be obtained from the boundary complex $\partial C_d^{\Delta}$
of the crosspolytope $C_d^{\Delta}$ with $2d$ vertices
by Jockusch's construction: We start with $\partial C_d^{\Delta}$
and compose a simplicial ball $B_{2d}$ 
as follows. Let the $(d-1)$-simplex\, $1\cdots d$\, belong to $B_{2d}$
and also all $d$-simplices\, $1\cdots k_1^I\cdots k_j^I\cdots d$, where
for\, $j=1,\ldots ,\frac{d-2}{2}$ the numbers $1\leq k_1<\ldots <k_j\leq d$\, are replaced by their images under
the involution\, $I=(1\,\,\,d+1)(2\,\,\,d+2)\cdots (d\,\,\,2d)$.
This collection of simplices $B_{2d}$ forms indeed a ball
(with boundary consisting of all $(d-2)$-faces
$1\cdots k_1^I\cdots\widehat{s}\cdots k_{(d-2)/2}^I\cdots d$\,
with vertex\, $s\in \{ 1,\ldots ,d\}$, $s\neq k_i$, deleted).
Moreover, $B_{2d}$ and $B^I_{2d}$ have the desired property
that 
\begin{itemize}
\item every $i$-face, $0\leq i\leq\lfloor \frac{d}{2}\rfloor -2$, of $\partial C_d^{\Delta}$
is contained in the boundaries of the two balls,
\item but no $(\lfloor\frac{d}{2}\rfloor -1)$-face of $\partial C_d^{\Delta}$
occurs as an interior face of the two balls.
\end{itemize}
If we remove the balls $B_{2d}$ and $B^I_{2d}$ 
from $\partial C_d^{\Delta}$ and sew in the new balls \linebreak
$(2d+1)*\partial B_{2d}$\, and\, $(2d+2)*\partial B^I_{2d}$, then
the resulting sphere is centrally symmetric and nearly neighborly.
In fact, it is isomorphic to $\partial H_{\,2d+2}^{\,d}$.

Besides the odd-dimensional polytopal spheres $\partial H_{\,2d+2}^{\,d}$, 
Bj\"orner, Paffenholz, Sj\"ostrand, and Ziegler \cite{BjoernerPaffenholzSjoestrandZiegler2003pre}
have recently constructed asymtotically many even-dimensional 
non-polytopal nearly neighborly centrally symmetric $(d-1)$-spheres
with $2d+2$ vertices that are Bier spheres.

\bigskip

Let us summarize the unsatisfactory present situation that we have
for centrally symmetric polytopes and spheres:
\begin{quote}
   \emph{Stanley \cite{Stanley1987} (and Novik \cite{Novik1999})
   proved a lower bound theorem for centrally symmetric polytopes, 
   but not for centrally symmetric spheres.}
\end{quote}
\begin{quote}
   \emph{Gr\"unbaum's centrally symmetric upper bound conjecture \cite{Gruenbaum1970a} 
   might well hold for spheres (but is wrong for polytopes).}
\end{quote}
\begin{quote}
   \emph{There are nearly neighborly centrally symmetric $d$-polytopes
   with $2d+2$ vertices (McMullen and Shephard~\cite{McMullenShephard1968})
   and nearly neighborly centrally symmetric $3$-spheres with $n=2k\geq 8$ vertices (Jockusch \cite{Jockusch1995}),
   but not much is known beyond these examples.}
\end{quote}
\begin{quote}
    \emph{According to Burton~\cite{Burton1991}, centrally symmetric $d$-polytopes 
    with sufficiently many vertices cannot be centrally $2$-neighborly.}
\end{quote}
In view of the result of Burton, presently not even a good guess for an upper bound
conjecture for centrally symmetric polytopes is available. Moreover,
we severly lack constructions that yield centrally symmetric polytopes
or spheres with many faces.

\subsection*{2\,\, Enumeration Results for Nearly Neighborly Spheres}

One approach to obtain nearly neighborly centrally 
symmetric spheres, at least on few vertices, 
is by computer enumeration. In \cite{BokowskiBremnerLutzMartin2003pre},
combinatorial $3$-mani\-folds are enumerated up to $10$ vertices.

\begin{thm} \cite{BokowskiBremnerLutzMartin2003pre}
There are exactly two non-isomorphic nearly neighborly centrally symmetric $3$-spheres
with $n=10$ vertices, the Gr\"unbaum sphere $G_{10}^{4}$ and
the Jockusch sphere $J^{\,3}_{10}$. 
\end{thm}

With the present enumeration techniques, an enumeration of \emph{all} nearly neighborly centrally symmetric 
$3$-spheres with $12$ vertices is already far out of reach. However, results for larger numbers of vertices 
can be achieved by restricting the enumeration to more symmetric triangulations.

In \cite{Lutz1999} we enumerated combinatorial $3$-manifolds
with a vertex-transitive automorphism group on up to $15$ vertices 
and found, besides $\partial C_4^{\Delta}$ and the Gr\"unbaum sphere $G_{10}^{4}$, 
two vertex-transitive nearly neighborly centrally symmetric $3$-spheres with
$12$ vertices and one with $14$ vertices. Apart from one example with
$12$-vertices, these spheres have a transitive cyclic automorphism group.
It therefore seemed promising to search for nearly neighborly
centrally symmetric spheres with a vertex-transitive 
cyclic (or dihedral) group action on more vertices and in higher dimensions $d$.

The standard dihedral and cyclic group action on the set $\{ 1,\dots,2k\}$,
with generators $a_{2k}=(123\ldots 2k)$\, and\, $b_{2k}=(1\,\,\, 2k)(2\,\,\, 2k\!-\!1)\ldots (k\,\,\,k\!+\!1)$\,
of\, $D_{2k}=\langle a_{2k},b_{2k}\rangle$\, and\, ${\mathbb Z}_{2k}=\langle a_{2k}\rangle$,
respectively, bring along a large number of small orbits of $(d+1)$-sets.
However, many of these orbits can be neglected if we are interested 
in centrally symmetric triangulations only:\linebreak 
We delete all orbits containing facets $F$ for which
$F\cap\,F^I\neq\emptyset$, with respect to the
involution $I=(12\cdots 2k)^k=(1\,\,k+1)\cdots (k\,\,2k)$,
in a preprocessing step 
before starting the enumeration program MANIFOLD\_VT \cite{Lutz_MANIFOLD_VT}. 
Every nearly neighborly centrally symmetric
example that we find we label with a unique symbol\,
$\mbox{}^d_{nn}\hspace{.3pt}n^{\,di/cy}_{\,z}$\, 
denoting the $z$-th isomorphism type of a \emph{n}early \emph{n}eighborly
centrally symmetric $d$-sphere listed for the 
dihedral/cyclic group action on $n=2k$ vertices.
For fixed $d$ and $n=2k$, we first process the dihedral
and then the cyclic action. The described search was carried out
in \cite[Ch.\ 4]{Lutz1999} for $3$-spheres with up to $16$ vertices
and has since then be extended to $22$ vertices.

\begin{table}[h]
\small\centering
\defaultaddspace=0.2em
\caption{Nearly neighborly centrally symmetric spheres with cyclic symmetry.}\label{tbl:enumeration}
\begin{tabular*}{\linewidth}{@{\extracolsep{\fill}}crrrrrrrrr@{}}
\\\toprule
 \addlinespace
$d\,\,\backslash\,\,n$ &  6 &  8 & 10 & 12 & 14 & 16 & 18 & 20 & 22 \\ \midrule
 \addlinespace
        2       &  1 &  0 &  0 &  0 &  0 &  0 &  0 &  0 &  0 \\
 \addlinespace
        3       & -- &  1 &  1 &  1 &  1 &  5 & 10 &  9 & 12 \\
 \addlinespace
        4       & -- & -- &  1 &  0 &  0 &  ? &  ? &  ? &  ? \\
 \addlinespace
        5       & -- & -- & -- &  1 &  2 &  3 &  ? &  ? &  ? \\
 \addlinespace
        6       & -- & -- & -- & -- &  1 &  0 &  ? &  ? &  ? \\
 \addlinespace
        7       & -- & -- & -- & -- & -- &  1 & 12 &  ? &  ? \\
 \addlinespace
 \addlinespace
\bottomrule
\end{tabular*}
\end{table}

\begin{thm}
There are nearly neighborly centrally symmetric $3$-spheres with a
vertex-transitive cyclic group action on $n=2k$ vertices for $4\leq k\leq 11$.
Moreover, there are nearly neighborly centrally symmetric $d$-spheres with a
vertex-transitive cyclic group action on $n=2k$ vertices for
$(d,n)=(5,14)$, $(5,16)$, $(7,18)$, but none for 
$(d,n)=(4,12)$, $(4,14)$, $(6,16)$.
(Table~\ref{tbl:enumeration} gives the respective numbers of spheres found by enumeration.)
\end{thm}

If $d=2$, then the boundaries of the tetrahedron, the octahedron,
and the icosahedron are the only vertex-transitive triangulations 
of the $2$-sphere $S^2$:
By Euler's formula, $f_0-f_1+f_2=2$, and double counting, $2f_1=3f_2$,
it follows that every triangulated $2$-sphere with $n$ vertices
has $f$-vector $f=(n,3n-6,2n-4)$. 
If the triangulation is 
vertex-transitive, then every vertex has the same number, say $q$,
of neighbors and is contained in exactly $q$ triangles. Double counting
yields $2f_1=nq$, or, equivalently, $(6-q)n=12$. The last equation
has three non-negative solutions $(n,q)=(4,3)$, $(6,4)$, and $(12,5)$.
The only possible examples corresponding to these values are
the boundaries of the tetrahedron, octahedron, and icosahedron. 
In particular, it follows that the boundary of the octahedron
is the only centrally symmetric $2$-sphere with a vertex-transitive cyclic group action.

\bigskip

\noindent
\textbf{Centrally Symmetric Cyclic Upper Bound Conjecture}
\emph{For all odd dimensions $d-1\geq 1$ and even\, $n=2k\geq 2d$,
there is a nearly neighborly centrally symmetric $(d-1)$-sphere 
with a vertex-trans\-itive cyclic group action on $n$ vertices.}

\bigskip

The conjecture is trivial for $d-1=1$ and clearly implies Gr\"unbaum's 
upper bound conjecture for centrally symmetric spheres
in odd, but also in even dimensions.
(The latter follows by suspending the respective
odd-dimensional examples.)

\begin{conj}
If $d$ is even, then the boundary complex of the $d$-dimensional crosspolytope
on $n=2d$ vertices is the only nearly neighborly centrally symmetric $d$-sphere
with a vertex-transitive cyclic group action.
\end{conj}

In Table~\ref{tbl:neighborly-spheres}, we list some of the spheres that 
we found by enumeration. The complete list of spheres is available
online at \cite{Lutz_PAGE}. If a sphere is centrally $l$-neighborly, 
i.e., if it has the $(l-1)$-skeleton of the corresponding cross-polytope, 
then we display the entry $f_l$ in italics (the entry
$n=f_0$ of the $f$-vector is listed separately in Column~2 of the table).
In Column~5 we list the respective orbit generators together with the 
corresponding orbit sizes as subscripts.

For some of the examples their full combinatorial automorphism group 
is larger than the dihedral or cyclic symmetry, indicated 
by the superscript $di$ or $cy$ in Table~\ref{tbl:neighborly-spheres}.
However, only few of the examples admit a dihedral symmetry.

\begin{landscape}

\small
\defaultaddspace=.1em

\setlength{\LTleft}{0pt}
\setlength{\LTright}{0pt}
\begin{longtable}{@{\extracolsep{\fill}}rrllll@{}}
\caption{\protect\parbox[t]{15cm}{Nearly neighborly centrally symmetric spheres with dihedral/cyclic group action.}}\label{tbl:neighborly-spheres}
\\\toprule
 $d$ & $n$ &   $f$-vector  & Type                              & List of orbits & Remarks \\ \midrule
\endfirsthead
\caption{\protect\parbox[t]{15cm}{Nearly neighborly centrally symmetric spheres (continued).}}
\\\toprule
 $d$ & $n$ &   $f$-vector  & Type                              & List of orbits & Remarks\\ \midrule
\endhead
\bottomrule
\endfoot
 \addlinespace
 \addlinespace
 \addlinespace
 \addlinespace
  2  &  6  &   (\emph{12},\emph{8})       & $\mbox{}^2_{nn}\hspace{.3pt}6^{\,di}_{\,1}$ & $123_{\,6}$\, $135_{\,2}$   & $\partial C_3^{\Delta}$, \\
     &     &                &&& \cite[$\mbox{}^2\hspace{.3pt}6^{\,11}_{\,1}$]{Lutz1999} \\
 \addlinespace
  3  &  8  &   (\emph{24},\emph{32},\emph{16})   & $\mbox{}^3_{nn}\hspace{.3pt}8^{\,di}_{\,1}$ & $1234_{\,8}$\, $1247_{\,8}$ & $\partial C_4^{\Delta}$, $C\!S^3_{8}$, \\
     &     &                &&& \cite[$\mbox{}^3\hspace{.3pt}8^{\,44}_{\,1}$]{Lutz1999} \\
 \addlinespace
     & 10  &   (\emph{40},60,30)   & $\mbox{}^3_{nn}\hspace{.3pt}10^{\,di}_{\,1}$& $1234_{\,10}$\, $1245_{\,10}$\, $1258_{\,10}$ & \cite[$\mbox{}^3\hspace{.3pt}10^{\,22}_{\,1}$]{Lutz1999} \\
 \addlinespace
     & 12  &   (\emph{60},96,48)   & $\mbox{}^3_{nn}\hspace{.3pt}12^{\,cy}_{\,1}$& $1234_{\,12}$\, $1246_{\,12}$\, $126\,11_{\,12}$\, $135\,10_{\,12}$ & $C\!S^3_{12}$, \cite[$\mbox{}^3\hspace{.3pt}12^{\,1}_{\,1}$]{Lutz1999} \\
 \addlinespace
     & 14  &  (\emph{84},140,70)   & $\mbox{}^3_{nn}\hspace{.3pt}14^{\,cy}_{\,1}$& $1234_{\,14}$\, $1245_{\,14}$\, $125\,10_{\,14}$\, $126\,10_{\,14}$\, $126\,12_{\,14}$ & \cite[$\mbox{}^3\hspace{.3pt}14^{\,1}_{\,1}$]{Lutz1999} \\
 \addlinespace
     & 16  & (\emph{112},192,96)   & $\mbox{}^3_{nn}\hspace{.3pt}16^{\,cy}_{\,1}$& $1234_{\,16}$\, $1246_{\,16}$\, $1268_{\,16}$\, $128\,15_{\,16}$\, $135\,14_{\,16}$\, $13\,10\,13_{\,16}$ & $C\!S^3_{16}$, \cite[$\mbox{}^3\hspace{.3pt}16^{\,1}_{\,56}$]{Lutz2003pre} \\
 \addlinespace
     &     &                & $\mbox{}^3_{nn}\hspace{.3pt}16^{\,cy}_{\,2}$& $1234_{\,16}$\, $1248_{\,16}$\, $1268_{\,16}$\, $126\,15_{\,16}$\, $1357_{\,16}$\, $138\,10_{\,16}$ & \cite[$\mbox{}^3\hspace{.3pt}16^{\,55}_{\,6}$]{Lutz2003pre} \\
 \addlinespace
     &     &                & $\mbox{}^3_{nn}\hspace{.3pt}16^{\,cy}_{\,3}$& $1234_{\,16}$\, $1248_{\,16}$\, $128\,15_{\,16}$\, $135\,12_{\,16}$\, $135\,14_{\,16}$\, $137\,14_{\,16}$ & \cite[$\mbox{}^3\hspace{.3pt}16^{\,1}_{\,58}$]{Lutz2003pre} \\
 \addlinespace
     &     &                & $\mbox{}^3_{nn}\hspace{.3pt}16^{\,cy}_{\,4}$& $1234_{\,16}$\, $124\,15_{\,16}$\, $1357_{\,16}$\, $136\,10_{\,16}$\, $137\,14_{\,16}$\, $13\,10\,13_{\,16}$ & \cite[$\mbox{}^3\hspace{.3pt}16^{\,1}_{\,63}$]{Lutz2003pre} \\
 \addlinespace
     &     &                & $\mbox{}^3_{nn}\hspace{.3pt}16^{\,cy}_{\,5}$& $1237_{\,16}$\, $1238_{\,16}$\, $126\,15_{\,16}$\, $128\,15_{\,16}$\, $1357_{\,16}$\, $13\,10\,13_{\,16}$ & \cite[$\mbox{}^3\hspace{.3pt}16^{\,55}_{\,9}$]{Lutz2003pre} \\
 \addlinespace
     & 18  & (\emph{144},252,126)   & $\mbox{}^3_{nn}\hspace{.3pt}18^{\,cy}_{\,1}$& $1234_{\,18}$\, $1245_{\,18}$\, $1256_{\,18}$\, $126\,12_{\,18}$\, $128\,12_{\,18}$\, $128\,15_{\,18}$\, $159\,13_{\,18}$ & \\
 \addlinespace
     &     &                & $\mbox{}^3_{nn}\hspace{.3pt}18^{\,cy}_{\,2}$& $1234_{\,18}$\, $1248_{\,18}$\, $128\,12_{\,18}$\, $12\,12\,17_{\,18}$\, $137\,14_{\,18}$\, $147\,11_{\,18}$\, $147\,14_{\,18}$ & \\
 \addlinespace
     &     &                & $\mbox{}^3_{nn}\hspace{.3pt}18^{\,cy}_{\,3}$& $1234_{\,18}$\, $1249_{\,18}$\, $125\,13_{\,18}$\, $125\,17_{\,18}$\, $129\,15_{\,18}$\, $12\,13\,15_{\,18}$\, $147\,12_{\,18}$ & \\
 \addlinespace
     &     &                & $\mbox{}^3_{nn}\hspace{.3pt}18^{\,cy}_{\,4}$& $1234_{\,18}$\, $1249_{\,18}$\, $126\,13_{\,18}$\, $126\,17_{\,18}$\, $129\,16_{\,18}$\, $12\,13\,16_{\,18}$\, $147\,12_{\,18}$ & \\
 \addlinespace
     &     &                & $\mbox{}^3_{nn}\hspace{.3pt}18^{\,cy}_{\,5}$& $1234_{\,18}$\, $124\,14_{\,18}$\, $126\,14_{\,18}$\, $126\,17_{\,18}$\, $138\,13_{\,18}$\, $147\,11_{\,18}$\, $147\,14_{\,18}$ & \\
 \addlinespace
     &     &                & $\mbox{}^3_{nn}\hspace{.3pt}18^{\,cy}_{\,6}$& $1234_{\,18}$\, $124\,15_{\,18}$\, $125\,15_{\,18}$\, $125\,17_{\,18}$\, $137\,14_{\,18}$\, $147\,14_{\,18}$\, $14\,11\,15_{\,18}$ & \\
 \addlinespace
     &     &                & $\mbox{}^3_{nn}\hspace{.3pt}18^{\,cy}_{\,7}$& $1235_{\,18}$\, $1236_{\,18}$\, $1249_{\,18}$\, $126\,14_{\,18}$\, $129\,14_{\,18}$\, $1358_{\,18}$\, $149\,15_{\,18}$ & \\
 \addlinespace
     &     &                & $\mbox{}^3_{nn}\hspace{.3pt}18^{\,cy}_{\,8}$& $1235_{\,18}$\, $1236_{\,18}$\, $124\,14_{\,18}$\, $1269_{\,18}$\, $129\,14_{\,18}$\, $1358_{\,18}$\, $138\,13_{\,18}$ & \\
 \addlinespace
     &     &                & $\mbox{}^3_{nn}\hspace{.3pt}18^{\,cy}_{\,9}$& $1236_{\,18}$\, $1237_{\,18}$\, $1247_{\,18}$\, $1248_{\,18}$\, $125\,12_{\,18}$\, $128\,12_{\,18}$\, $147\,11_{\,18}$ & \\
 \addlinespace
     &     &                & $\mbox{}^3_{nn}\hspace{.3pt}18^{\,cy}_{\,10}$& $1236_{\,18}$\, $1237_{\,18}$\, $125\,17_{\,18}$\, $127\,12_{\,18}$\, $12\,12\,17_{\,18}$\, $1369_{\,18}$\, $139\,14_{\,18}$ & \\
 \addlinespace
\newpage
 \addlinespace
 \addlinespace
 \addlinespace
 \addlinespace
     & 20  & (\emph{180},320,160)   & 
$\mbox{}^3_{nn}\hspace{.3pt}20^{\,cy}_{\,1}$ &  & \cite{Lutz_PAGE}\\ 
     &     &                & -- $\mbox{}^3_{nn}\hspace{.3pt}20^{\,cy}_{\,9}$ &  & \\
 \addlinespace
     & 22  & (\emph{220},396,198)   & 

$\mbox{}^3_{nn}\hspace{.3pt}22^{\,cy}_{\,1}$ &  & \cite{Lutz_PAGE}\\ 
     &     &                & -- $\mbox{}^3_{nn}\hspace{.3pt}22^{\,cy}_{\,12}$ &  & \\
 \addlinespace
  4  & 10  & (\emph{40},\emph{80},\emph{80},\emph{32})  & $\mbox{}^4_{nn}\hspace{.3pt}10^{\,di}_{\,1}$& $12345_{\,10}$\, $12359_{\,10}$\, $12458_{\,10}$\, $13579_{\,2}$ & $\partial C_5^{\Delta}$, \\
     &     &                &&& \cite[$\mbox{}^4\hspace{.3pt}10^{\,39}_{\,1}$]{Lutz1999} \\
 \addlinespace
  5  & 12  & (\emph{60},\emph{160},\emph{240}, &$\mbox{}^5_{nn}\hspace{.3pt}12^{\,di}_{\,1}$& $123456_{\,12}$\, $12346\,11_{\,12}$\, $12356\,10_{\,24}$ & $\partial C_6^{\Delta}$, \\ 
     &     & \emph{192},\emph{64})        &                                   & $12469\,11_{\,12}$\, $12569\,10_{\,4}$ & \cite[$\mbox{}^5\hspace{.3pt}12^{\,293}_{\,1}$]{Lutz1999} \\
 \addlinespace
     & 14  & (\emph{84},\emph{280},490, &$\mbox{}^5_{nn}\hspace{.3pt}14^{\,di}_{\,1}$& $123456_{\,14}$\, $123467_{\,28}$\, $1234712_{\,14}$\, $12367\,12_{\,28}$ & \cite[$\mbox{}^5\hspace{.3pt}14^{\,49}_{\,1}$]{Lutz1999} \\ 
     &     & 420,140)       &                                   & $12457\,10_{\,28}$\, $1247\,10\,13_{\,14}$\, $1256\,10\,11_{\,14}$ & \\
 \addlinespace
     &     &                &$\mbox{}^5_{nn}\hspace{.3pt}14^{\,di}_{\,2}$& $123456_{\,14}$\, $12346\,12_{\,28}$\, $1234712_{\,14}$\, $12356\,11_{\,28}$ & \cite[$\mbox{}^5\hspace{.3pt}14^{\,7}_{\,1}$]{Lutz1999} \\ 
     &     &                &                                   & $12467\,10_{\,28}$\, $1247\,10\,13_{\,14}$\, $1256\,10\,11_{\,14}$ & \\
 \addlinespace
     & 16  & (\emph{112},\emph{448},864, &$\mbox{}^5_{nn}\hspace{.3pt}16^{\,cy}_{\,1}$& $123456_{\,16}$\, $123467_{\,16}$\, $123478_{\,16}$\, $12348\,13_{\,16}$ & \\ 
     &     & 768,256)       &                                   & $1234\,13\,15_{\,16}$\, $12378\,12_{\,16}$\, $1238\,12\,13_{\,16}$\, $123\,12\,13\,15_{\,16}$ & \\ 
     &     &                &                                   & $123\,12\,14\,15_{\,16}$\, $12468\,11_{\,16}$\, $1246\,11\,15_{\,16}$\, $1248\,11\,13_{\,16}$ & \\ 
     &     &                &                                   & $124\,11\,13\,15_{\,16}$\, $1267\,11\,12_{\,16}$\, $1268\,11\,13_{\,16}$\, $1358\,10\,14_{\,16}$ & \\
 \addlinespace
     &     &                &$\mbox{}^5_{nn}\hspace{.3pt}16^{\,cy}_{\,2}$& $123456_{\,16}$\, $123467_{\,16}$\, $123478_{\,16}$\, $12348\,13_{\,16}$ & \\ 
     &     &                &                                   & $1234\,13\,15_{\,16}$\, $12378\,12_{\,16}$\, $1238\,12\,13_{\,16}$\, $123\,12\,13\,15_{\,16}$ & \\ 
     &     &                &                                   & $123\,12\,14\,15_{\,16}$\, $12468\,15_{\,16}$\, $1248\,13\,15_{\,16}$\, $1267\,11\,12_{\,16}$ & \\ 
     &     &                &                                   & $1268\,11\,13_{\,16}$\, $1268\,11\,15_{\,16}$\, $128\,11\,13\,15_{\,16}$\, $1358\,12\,14_{\,16}$ & \\
 \addlinespace 
     &     &                &$\mbox{}^5_{nn}\hspace{.3pt}16^{\,cy}_{\,3}$& $123456_{\,16}$\, $123467_{\,16}$\, $123478_{\,16}$\, $12348\,13_{\,16}$ & \\ 
     &     &                &                                   & $1234\,13\,15_{\,16}$\, $12378\,12_{\,16}$\, $1238\,12\,15_{\,16}$\, $1238\,13\,15_{\,16}$ & \\ 
     &     &                &                                   & $123\,12\,14\,15_{\,16}$\, $12467\,13_{\,16}$\, $1246\,11\,13_{\,16}$\, $1246\,11\,15_{\,16}$ & \\ 
     &     &                &                                   & $12478\,13_{\,16}$\, $124\,11\,13\,15_{\,16}$\, $1267\,11\,12_{\,16}$\, $1357\,10\,12_{\,16}$ & \\
 \addlinespace
\newpage
 \addlinespace
 \addlinespace
 \addlinespace
 \addlinespace
  6  & 14  &(\emph{84},\emph{280},\emph{560},&$\mbox{}^6_{nn}\hspace{.3pt}14^{\,di}_{\,1}$& $1234567_{\,14}$\, $123457\,13_{\,14}$\, $123467\,12_{\,28}$ & $\partial C_7^{\Delta}$, \\ 
     &     &\emph{672},\emph{448},\emph{128})&                                   & $123567\,11_{\,14}$\, $12357\,11\,13_{\,14}$\, $12367\,11\,12_{\,14}$ & \cite[$\mbox{}^6\hspace{.3pt}14^{\,57}_{\,1}$]{Lutz1999} \\ 
     &     &                &                                   & $12457\,10\,13_{\,14}$\, $12467\,10\,12_{\,14}$\, $13579\,11\,13_{\,2}$ & \\
 \addlinespace
  7  & 16  &(\emph{112},\emph{448},\emph{1120},& $\mbox{}^7_{nn}\hspace{.3pt}16^{\,di}_{\,1}$ & $12345678_{\,16}$\, $1234568\,15_{\,16}$\, $1234578\,14_{\,32}$\, $1234678\,13_{\,32}$ & $\partial C_8^{\Delta}$ \\
     &     &\emph{1792},\emph{1792},&                                   & $123468\,13\,15_{\,16}$\, $123478\,13\,14_{\,16}$\, $123568\,12\,15_{\,32}$\, $123578\,12\,14_{\,32}$ & \\
     &     &\emph{1024},\emph{256})&                                   & $123678\,12\,13_{\,16}$\, $124578\,11\,14_{\,16}$\, $12468\,11\,13\,15_{\,16}$\, $12478\,11\,13\,14_{\,16}$ & \\
 \addlinespace
     & 18  &(\emph{144},\emph{672},\emph{2016}, & 
$\mbox{}^7_{nn}\hspace{.3pt}18^{\,cy}_{\,1}$ &  & \cite{Lutz_PAGE}\\ 
     &     & 3780,4200,     & -- $\mbox{}^7_{nn}\hspace{.3pt}18^{\,cy}_{\,10}$ &  & \\
     &     & 2520,630)      & &  & \\
 \addlinespace
     &     &                & $\mbox{}^7_{nn}\hspace{.3pt}18^{\,di}_{\,1}$ & $12345678_{\,18}$\, $12345689_{\,36}$\, $1234569\,16_{\,18}$\, $1234589\,16_{\,36}$ & \\
     &     &                &                                   & $12346789_{\,18}$\, $1234679\,14_{\,36}$\, $123467\,14\,17_{\,36}$\, $123469\,14\,17_{\,18}$ & \\
     &     &                &                                   & $1234789\,14_{\,36}$\, $123478\,14\,15_{\,36}$\, $123489\,14\,15_{\,18}$\, $1235689\,13_{\,36}$ & \\
     &     &                &                                   & $123569\,13\,16_{\,36}$\, $123589\,13\,16_{\,36}$\, $123678\,13\,14_{\,18}$\, $123689\,13\,14_{\,36}$ & \\
     &     &                &                                   & $124578\,12\,15_{\,18}$\, $124579\,12\,15_{\,36}$\, $124589\,15\,16_{\,18}$\, $124679\,12\,17_{\,36}$ & \\
     &     &                &                                   & $12479\,12\,14\,15_{\,18}$\, $12479\,12\,14\,17_{\,18}$\, $12569\,12\,13\,16_{\,18}$ & \\
 \addlinespace
     &     &                & $\mbox{}^7_{nn}\hspace{.3pt}18^{\,di}_{\,2}$ & $12345678_{\,18}$\, $12345689_{\,36}$\, $1234569\,16_{\,18}$\, $1234589\,16_{\,36}$ & \\
     &     &                &                                   & $12346789_{\,18}$\, $1234679\,14_{\,36}$\, $123467\,14\,17_{\,36}$\, $123469\,14\,17_{\,18}$ & \\
     &     &                &                                   & $1234789\,14_{\,36}$\, $123478\,14\,15_{\,36}$\, $123489\,14\,15_{\,18}$\, $1235689\,16_{\,36}$ & \\
     &     &                &                                   & $123568\,13\,16_{\,36}$\, $123678\,13\,14_{\,18}$\, $123689\,13\,14_{\,36}$\, $123689\,13\,16_{\,36}$ & \\
     &     &                &                                   & $124578\,12\,15_{\,18}$\, $124589\,12\,15_{\,36}$\, $124589\,15\,16_{\,18}$\, $124679\,12\,17_{\,36}$ & \\
     &     &                &                                   & $12479\,12\,14\,15_{\,18}$\, $12479\,12\,14\,17_{\,18}$\, $12569\,12\,13\,16_{\,18}$ & \\
 \addlinespace
\end{longtable}

\end{landscape}


\subsection*{3\,\, A Transitive Series of Nearly Neighborly Spheres}

In this section, we prove the centrally symmetric cyclic upper bound conjecture
for $d=3$ for all numbers $n=4m\geq 8$ of vertices.

\begin{thm}\label{thm:4m_series}
There is an infinite series of nearly neighborly centrally symmetric 
$3$-spheres $C\!S^3_{4m}$ with a transitive cyclic group action on $4m$
vertices \linebreak
for $m\geq 2$. 
\end{thm}
\emph{Proof.}
Let the permutation $g=(1,2,\dots,4m)$ be the generator of the standard transitive cyclic group
action on the vertex set $\{ 1,2,\dots,4m\}$.
We define a series of $3$-dimensional simplicial complexes $C\!S^3_{4m}$ 
in terms of the orbit generators of Table~\ref{tbl:seriesT34m}:
Let every orbit generator $ijkl_{\,4m}$, with the orbit-size as index, 
contribute an orbit of $4m$ tetrahedral facets 
$ijkl$, $(i+1)(j+1)(k+1)(l+1)$, $\dots$, $(i+4m)(j+4m)(k+4m)(l+4m)$
to the simplicial complex $C\!S^3_{4m}$, where the vertex-labels
are to be taken modulo $4m$.

\bigskip

\begin{table}[h]
\small\centering
\defaultaddspace=0.2em
\caption{The series $C\!S^3_{4m}$.}\label{tbl:seriesT34m}
\begin{tabular}{@{}c@{\hspace{2mm}}l@{\hspace{2mm}}l@{\hspace{2mm}}l@{\hspace{2mm}}l@{}}
\\\toprule
 \addlinespace
 Sphere &\multicolumn{4}{@{}l}{List of Orbits} \\ \midrule
 \addlinespace
 $C\!S^3_{8}$ & $1234_{\,8}$  & $1247_{\,8}$ \\
 \addlinespace
 $C\!S^3_{12}$ & $1234_{\,12}$ & $1249_{\,12}$      & $129\,11_{\,12}$        & $1358_{\,12}$ \\
 \addlinespace
 $C\!S^3_{16}$ & $1234_{\,16}$ & $124\,11_{\,16}$   & $12\,11\,13_{\,16}$     & $1358_{\,16}$ \\
      &               &                             & $12\,13\,15_{\,16}$     & $137\,10_{\,16}$ \\
 \addlinespace
$\cdots$ \\
 \addlinespace
 $C\!S^3_{4m}$ & $1234_{\,4m}$ & $124(2m+3)_{\,4m}$ & $12(2m+3)(2m+5)_{\,4m}$ & $1358_{\,4m}$ \\
      &               &                             & $12(2m+5)(2m+7)_{\,4m}$ & $137\,10_{\,4m}$ \\
      &               &                             & $\cdots$                & $\cdots$ \\
      &               &                             & $12(4m-3)(4m-1)_{\,4m}$ & $13(2m-1)(2m+2)_{\,4m}$ \\
 \addlinespace
\bottomrule
\end{tabular}
\end{table}

\bigskip

By construction, $C\!S^3_{4m}$ is invariant under the standard vertex-transive
cyclic symmetry, in particular, it is invariant under the involution\,
$I:=$ \linebreak $(1,2,\dots,4m)^{2m}=(1,2m+1)(2,2m+2)\dots(2m,4m)$.
No (non-empty) face of $C\!S^3_{4m}$ is fixed under~$I$, which easily can be
verified by inspecting the defining orbits of $C\!S^3_{4m}$.
Hence, $C\!S^3_{4m}$ is a centrally symmetric $3$-dimensional simplicial complex.

In the following, we will prove that $C\!S^3_{4m}$ is a $3$-sphere 
by showing that $C\!S^3_{4m}$ is a $3$-manifold of Heegaard genus one 
with a Heegaard diagram that has one crossing (cf.\ \cite{KuehnelLutz2003pre}, \cite[Sec.\ 63]{SeifertThrelfall1980}). 
Moreover, we will see that $C\!S^3_{4m}$ is nearly neighborly.

\begin{figure}[h]
\begin{center}
\footnotesize
\psfrag{1}{1}
\psfrag{2}{2}
\psfrag{3}{3}
\psfrag{4}{4}
\psfrag{5}{5}
\psfrag{6}{6}
\psfrag{7}{7}
\psfrag{8}{8}
\includegraphics[width=.425\linewidth]{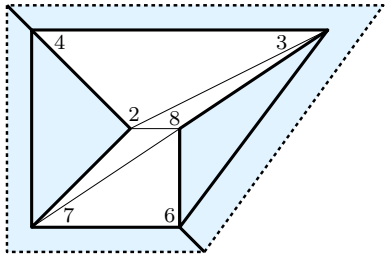}
\end{center}
\caption{The link of vertex 1 in $C\!S^3_8$.}
\label{fig:link8}
\end{figure}

\begin{figure}[h]
\begin{center}
\footnotesize
\psfrag{1}{1}
\psfrag{2}{2}
\psfrag{3}{3}
\psfrag{4}{4}
\psfrag{5}{5}
\psfrag{6}{6}
\psfrag{7}{7}
\psfrag{8}{8}
\psfrag{9}{9}
\psfrag{10}{10}
\psfrag{11}{11}
\psfrag{12}{12}
\psfrag{I}{\normalsize I}
\psfrag{II}{\normalsize II}
\psfrag{III}{\normalsize III}
\psfrag{IV}{\normalsize IV}
\includegraphics[width=.55\linewidth]{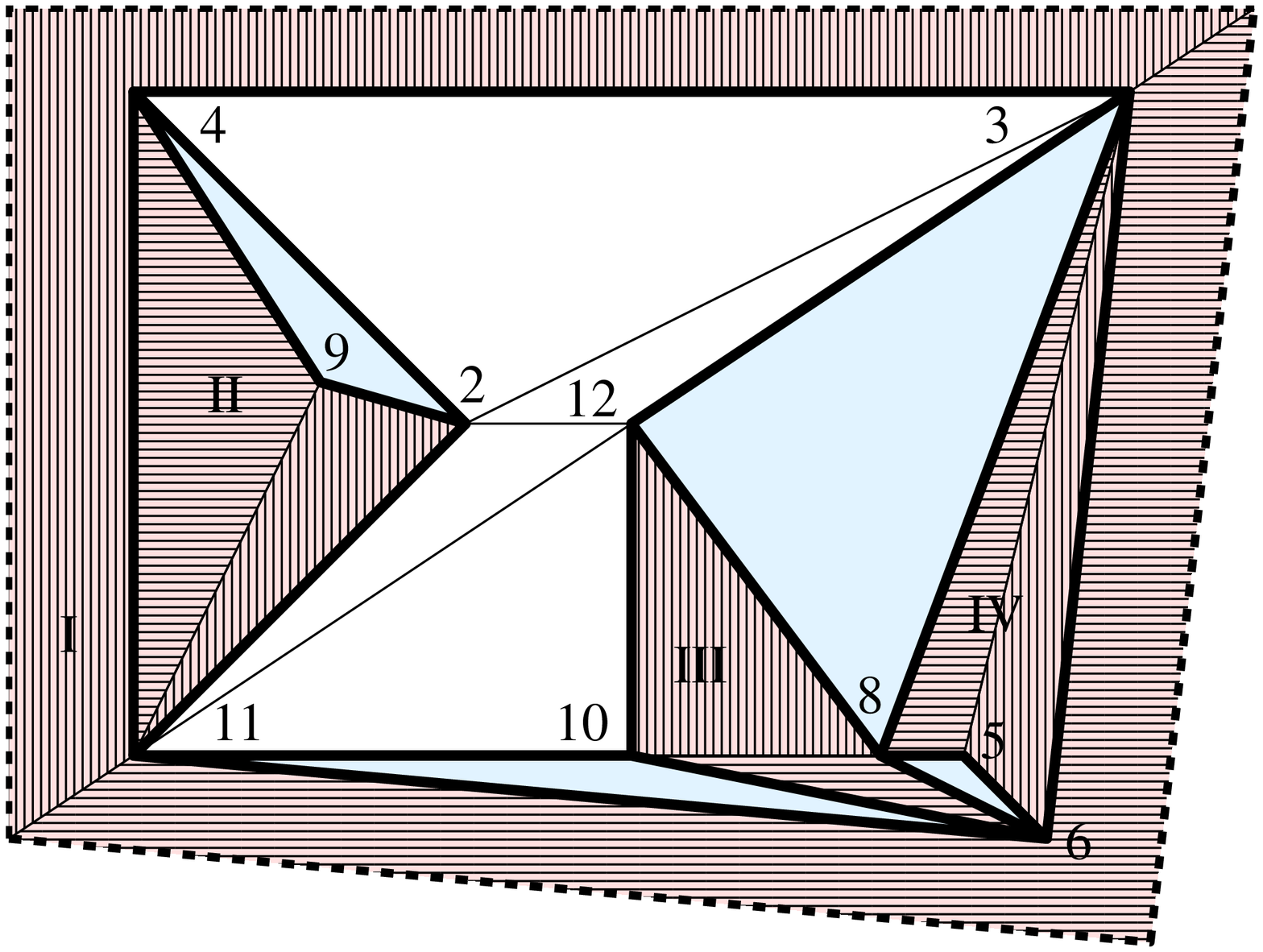}
\end{center}
\caption{The link of vertex 1 in $C\!S^3_{12}$.}
\label{fig:link12}
\end{figure}

\begin{figure}
\begin{center}
\footnotesize
\psfrag{1}{1}
\psfrag{2}{2}
\psfrag{3}{3}
\psfrag{4}{4}
\psfrag{5}{5}
\psfrag{6}{6}
\psfrag{7}{7}
\psfrag{8}{8}
\psfrag{9}{9}
\psfrag{10}{10}
\psfrag{11}{11}
\psfrag{12}{12}
\psfrag{13}{13}
\psfrag{14}{14}
\psfrag{15}{15}
\psfrag{16}{16}
\psfrag{I}{\normalsize I}
\psfrag{II}{\normalsize II}
\psfrag{III}{\normalsize III}
\psfrag{IV}{\normalsize IV}
\includegraphics[width=.675\linewidth]{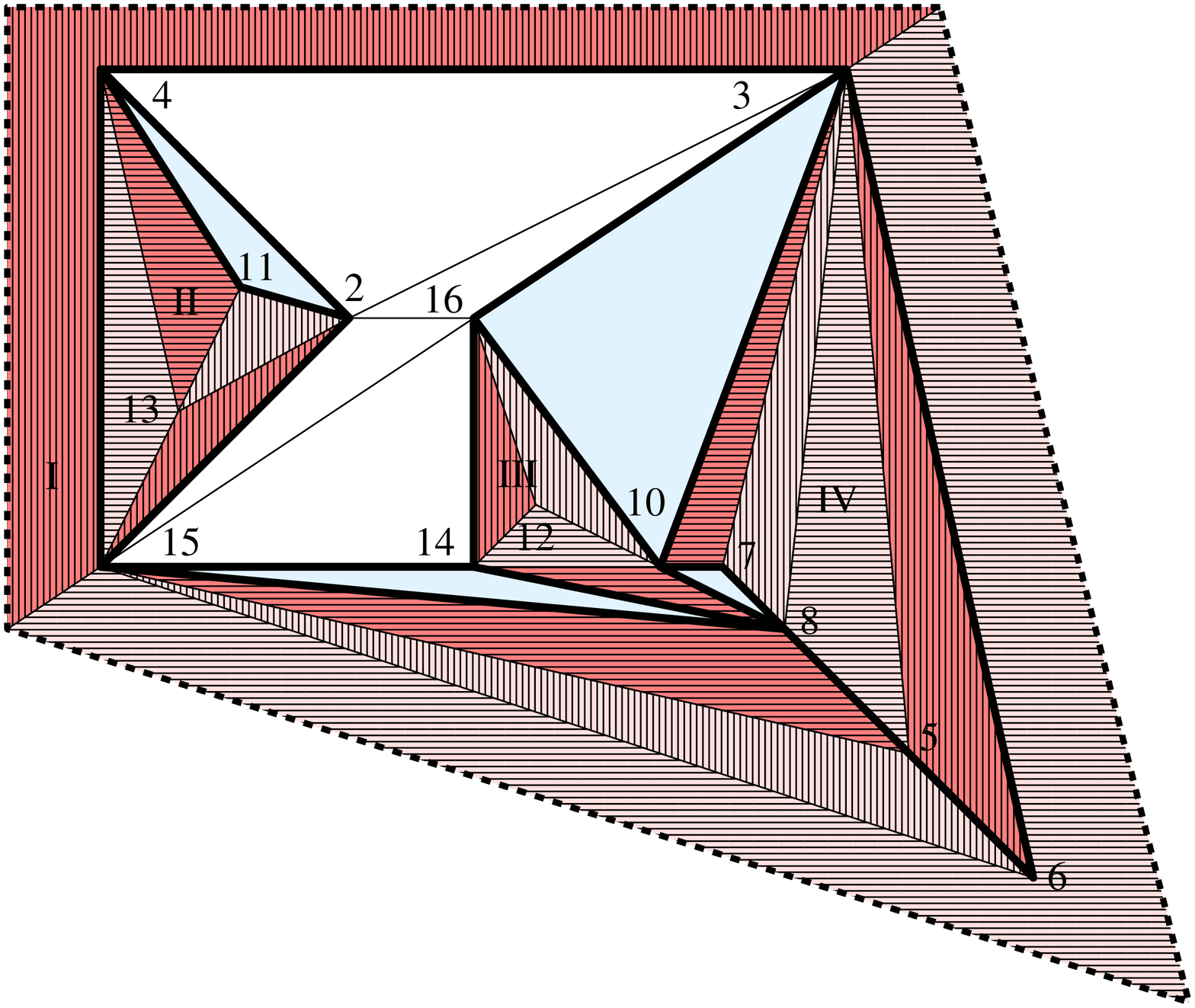}
\end{center}
\caption{The link of vertex 1 in $C\!S^3_{16}$.}
\label{fig:link16}
\end{figure}

\begin{figure}
\begin{center}
\footnotesize
\psfrag{1}{1}
\psfrag{2}{2}
\psfrag{3}{3}
\psfrag{4}{4}
\psfrag{5}{5}
\psfrag{6}{6}
\psfrag{7}{7}
\psfrag{8}{8}
\psfrag{9}{9}
\psfrag{10}{10}
\psfrag{11}{11}
\psfrag{12}{12}
\psfrag{13}{13}
\psfrag{14}{14}
\psfrag{15}{15}
\psfrag{16}{16}
\psfrag{17}{17}
\psfrag{18}{18}
\psfrag{19}{19}
\psfrag{20}{20}
\psfrag{I}{\normalsize I}
\psfrag{II}{\normalsize II}
\psfrag{III}{\normalsize III}
\psfrag{IV}{\normalsize IV}
\includegraphics[width=.8\linewidth]{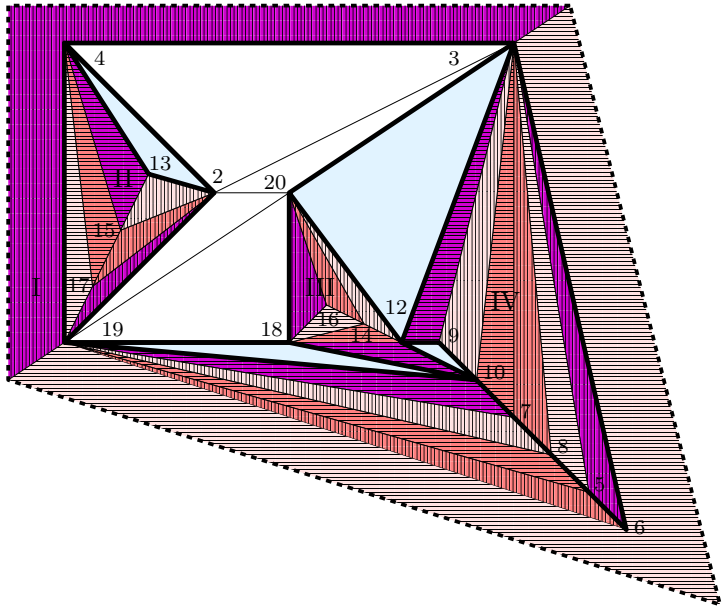}
\end{center}
\caption{The link of vertex 1 in $C\!S^3_{20}$.}
\label{fig:link20}
\end{figure}

In order to verify that $C\!S^3_{4m}$ is a $3$-manifold, we need to show
that the link of every of its vertices is a triangulated $2$-sphere.
Since $C\!S^3_{4m}$ is vertex-transitive, it suffices to analyze the link
of vertex~$1$. The vertex-links of vertex~$1$ in the complexes
$C\!S^3_8$, $C\!S^3_{12}$, $C\!S^3_{16}$, and $C\!S^3_{20}$
are depicted in the Figures~\ref{fig:link8}, \ref{fig:link12}, \ref{fig:link16},
and~\ref{fig:link20}, respectively. The complex $C\!S^3_{4m}$ consists
of $2m-2$ orbits that contribute four triangles each to the link
of vertex 1. The orbits can be grouped into four different types:
The basic orbits $1234_{\,4m}$ (contributing white triangles) and $124(2m+3)_{\,4m}$
(contributing shaded triangles) in the columns 2 and 3 of Table~\ref{tbl:seriesT34m}
and the two series of orbits in the columns 4 and 5 of Table~\ref{tbl:seriesT34m} 
(contributing triangles with vertical and horizontal stripes, respectively).
The striped triangles form four different regions I--IV of $2m-4$
triangles each, half of them vertically and half of them horizontally
striped, respectively. Topologically, each of the four regions is a disc, but displays
a different kind of ``cristallographic growth'' when we increase $m$.
For example, region II consists of the $m-2$ vertically striped triangles
$2(2m+3)(2m+5)$, $2(2m+5)(2m+7)$, $\dots$, $2(4m-3)(4m-1)$
and of the $m-2$ horizontally striped triangles 
$4(2m+3)(2m+5)$, $4(2m+5)(2m+7)$, $\dots$, $4(4m-3)(4m-1)$.
It is easy to check that the four regions I--IV together with 
the four white triangles and the four shaded triangles form a $2$-sphere.
Hence, $C\!S^3_{4m}$ is a $3$-manifold.


\begin{figure}[h]
\begin{center}
\footnotesize
\psfrag{1}{1}
\psfrag{2}{2}
\psfrag{3}{3}
\psfrag{4}{4}
\psfrag{5}{5}
\psfrag{6}{6}
\psfrag{7}{7}
\psfrag{8}{8}
\includegraphics[height=.23\linewidth]{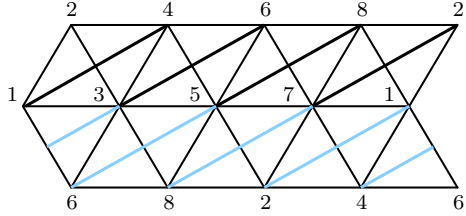}
\end{center}
\caption{The solid $3$-torus $T^3_8$.}
\label{fig:torusT_3_8}
\end{figure}

\begin{figure}[h]
\begin{center}
\footnotesize
\psfrag{1}{1}
\psfrag{2}{2}
\psfrag{3}{3}
\psfrag{4}{4}
\psfrag{5}{5}
\psfrag{6}{6}
\psfrag{7}{7}
\psfrag{8}{8}
\psfrag{9}{9}
\psfrag{10}{10}
\psfrag{11}{11}
\psfrag{12}{12}
\includegraphics[height=.23\linewidth]{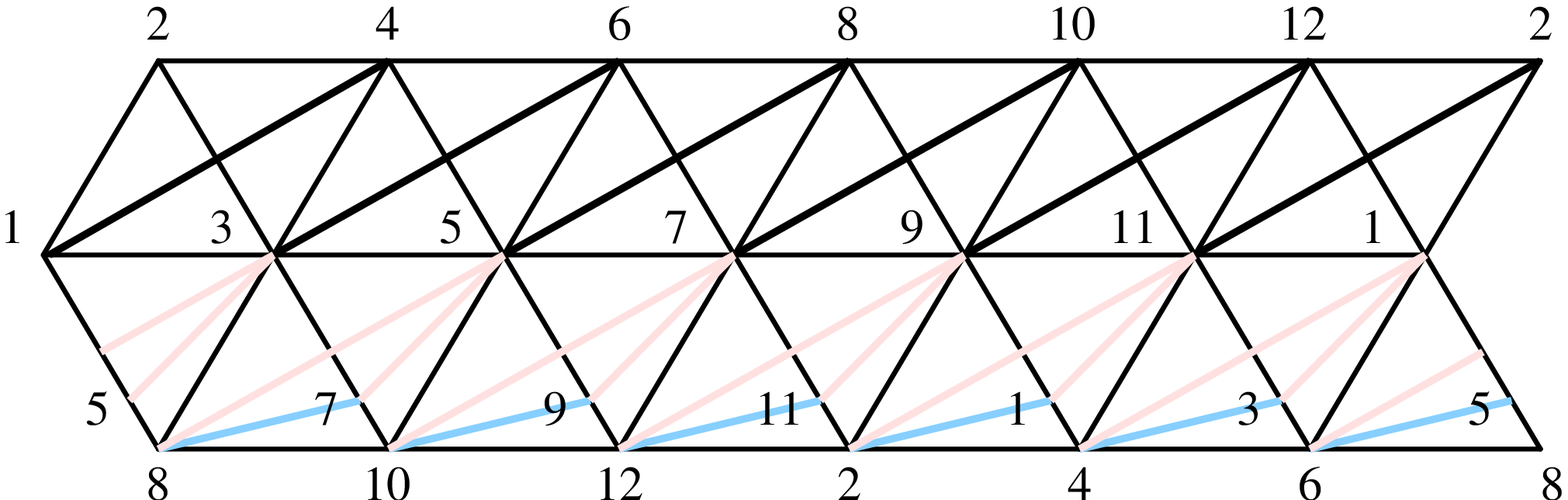}
\end{center}
\caption{The solid $3$-torus $T^3_{12}$.}
\label{fig:torusT_3_12}
\end{figure}

\begin{figure}[h]
\begin{center}
\footnotesize
\psfrag{1}{1}
\psfrag{2}{2}
\psfrag{3}{3}
\psfrag{4}{4}
\psfrag{5}{5}
\psfrag{6}{6}
\psfrag{7}{7}
\psfrag{8}{8}
\psfrag{9}{9}
\psfrag{10}{10}
\psfrag{11}{11}
\psfrag{12}{12}
\psfrag{13}{13}
\psfrag{14}{14}
\psfrag{15}{15}
\psfrag{16}{16}
\includegraphics[height=.23\linewidth]{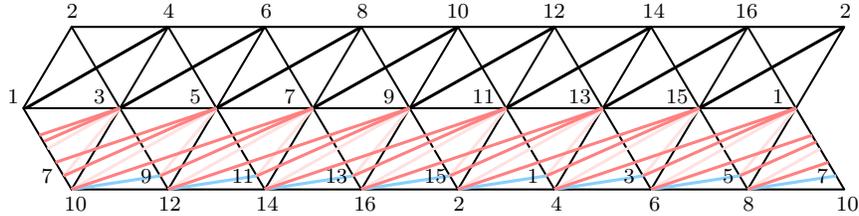}
\end{center}
\caption{The solid $3$-torus $T^3_{16}$.}
\label{fig:torusT_3_16}
\end{figure}

The triangulated $3$-manifold $C\!S^3_{4m}$ contains
as a $2$-dimensional subcomplex a vertex-transitive 
$2$-torus $T^2_{4m}$ with orbit generators $123_{\,4m}$ and $13(2m+2)_{\,4m}$.
We will show that this triangulated $2$-torus $T^2_{4m}$ splits $C\!S^3_{4m}$
into two parts, $T^3_{4m}$ and $(T^3_{4m})^g$, each of which 
is a triangulated solid $3$-torus and is mapped onto the 
other side by the glide reflection $g=(1,2,\dots,4m)$ 
of the $2$-torus $T^2_{4m}$. 
\begin{figure}[h]
\begin{center}
\footnotesize
\psfrag{1}{1}
\psfrag{2}{2}
\psfrag{3}{3}
\psfrag{4}{4}
\psfrag{5}{5}
\psfrag{6}{6}
\psfrag{7}{7}
\psfrag{8}{8}
\includegraphics[height=.24\linewidth]{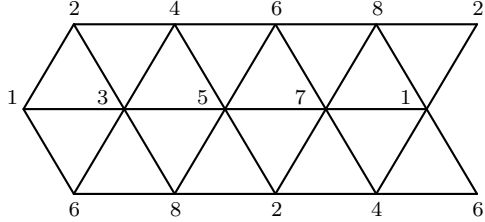}
\end{center}
\caption{The $2$-torus $T^2_8$.}
\label{fig:torusT_2_8}
\end{figure}
The $2$-torus $T^2_8$
is depicted in Figure~\ref{fig:torusT_2_8} with 
the orbits $123_8$ and $136_8$ forming the upper eight
and the lower eight triangles, respectively.
In Figures~\ref{fig:torusT_3_8}, \ref{fig:torusT_3_12}, and \ref{fig:torusT_3_16},
the tori $T^2_8$, $T^2_{12}$, and $T^2_{16}$ form the respective base grids.
In Figure~\ref{fig:torusT_3_8}, we glue ``on top'' of the upper eight triangles of $T^2_8$ 
every second tetrahedron of the orbit $1234_8$, i.e., 
the tetrahedra $1234$, $3456$, $5678$, and $1278$, 
as well as ``on top'' of the lower eight triangles of $T^2_8$  
every second tetrahedron of the orbit $1247_{\,8}$,
i.e., the tetrahedra $1247$, $1346$, $3568$, and $2578$. 
From the figure we see that every ``top'' triangular face
of one of the upper four tetrahedra appears also as a ``top'' triangular face
of one of the lower four tetrahedra. Hence, the tetrahedra of the upper half fit together
with the tetrahedra of the lower half to form a solid $3$-torus~$T^3_8$ whose boundary 
is, as the ``back side'', the torus $T^2_8$. 
In general, we also glue ``on top'' of the upper $4m$ triangles of $T^2_{4m}$ 
every second tetrahedron of the basic orbit $1234_{\,4m}$. ``On top'' of the
lower $4m$ triangles of $T^2_{4m}$, however, we first glue every
second tetrahedron of the basic orbit $124(2m+3)_{\,4m}$ and then
every second tetrahedron of the orbits alternatingly from the columns 5 and 4 
of Table~\ref{tbl:seriesT34m}. Upon completion, 
the ``top'' triangles of the upper part fit together with the ``top''
triangles of the lower part to form a
solid $3$-torus $T^3_{4m}$. Since $T^3_{4m}$ contains
every second tetrahedron of the orbits of $C\!S^3_{4m}$,
its image $(T^3_{4m})^g$ under the cyclic shift $g=(1,2,\dots,4m)$ 
has as its facets precisely the remaining tetrahedra of $C\!S^3_{4m}$
and, hence, is again a solid $3$-torus. Thus we have established that $C\!S^3_{4m}$
has a Heegaard splitting of genus one into the two solid tori
$T^3_{4m}$ and $(T^3_{4m})^g$.

The Heegaard diagram of $C\!S^3_{4m}$ consists of the
middle torus $T^2_{4m}$ together with a meridian circle $c$ of $T^3_{4m}$
and a meridian 
$c'$ of $(T^3_{4m})^g$. As meridian
of $T^3_{4m}$ we take 
$c:=(2m+1)(2m+3),\dots,(4m-3)(4m-1),(4m-1)(4m),(4m)(2m+1)$
on $T^2_{4m}$. Its image 
$c':=c^g=(2m+2)(2m+4),\dots,(4m-2)(4m),(4m)1$, $1(2m+2)$ 
under the glide reflection $g$ is a meridian of $(T^3_{4m})^g$
and intersects $c$ in the one 
crossing point $4m$. Since a $3$-manifold
$M$ is a $3$-sphere if it has a 
genus one 
Heegaard diagram with one crossing point, we are done.

\pagebreak

\mbox{}

\vspace{-.5mm}
\vspace{-.5mm}


It remains to show that the centrally symmetric $3$-sphere
$C\!S^3_{4m}$ is nearly neighborly.
Since the $f$-vector $(f_0,f_1,f_2,f_3)$ of 
a $3$-manifold is already determined by the number of vertices
$f_0$ and the number of facets $f_3$ via Euler's formula
$f_0-f_1+f_2-f_3=0$ and the Dehn-Sommerville equation
$f_2=2f_3$, it follows directly from the number and sizes of the defining orbits
that $C\!S^{\,3}_{\,4m}$ has $f$-vector $(4m,8m^2-4m,16m^2-16m,8m^2-8m)$.
Since $8m^2-4m=\binom{4m}{2}-2m$, the centrally symmetric 
$3$-sphere $C\!S^{\,3}_{\,4m}$ has the $1$-skeleton 
of the corresponding cross-polytope $C_{4m}^{\Delta}$ on $4m$ vertices
and, therefore, is nearly neighborly.\\ \mbox{} \hfill $\Box$

\begin{cor}
The nearly neighborly centrally symmetric $3$-spheres $C\!S^{\,3}_{\,4m}$
are not obtainable by Jockusch's construction for $m\geq 3$.
\end{cor}

\pagebreak

\emph{Proof.} The $3$-balls $B_{2k}$ in Jockusch's construction
are choosen such that they contain all vertices of $J^{\,3}_{2k}$,
but not the star of any edge of $J^{\,3}_{2k}$.
In particular, the boundary $2$-spheres $\partial B_{2k}$ are stacked spheres
and occur as the link of the vertices $2k+1$ in $J^{\,3}_{2k+2}$.
On the contrary, the vertex-links in the spheres $C\!S^{\,3}_{\,4m}$
are not stacked. \hfill $\Box$

\bigskip

Although the proof of correctness for the examples of Theorem~\ref{thm:4m_series} 
is rather straight forward, it is, in general, not at all obvious 
how we can find or construct \emph{series of vertex-transitive 
triangulations} of spheres or of other manifolds.
In the case of the series $C\!S^{\,3}_{\,4m}$ the generating orbits 
were discovered by examining the examples of Table~\ref{tbl:neighborly-spheres},
but all attempts failed so far to extend the series to or to find
alternative series on $4m+2$ vertices for $m\geq 2$.

Most surprising, however, is that we presently know of \emph{merely five
basic infinite series} of vertex-transitive triangulations of spheres:
\begin{itemize}
\item the boundary complexes of even-dimensional cyclic polytopes $C_d(n)$,
\item the boundary complexes of bicyclic $4$-polytopes $BiC(p,q;n)$
      of Smilansky \cite{Smilansky1990} for appropriate parameters $p$, $q$, and $n$  
      (cf.\ also \cite{BokowskiSchuchert1995} and \cite{Schuchert1995}),
\item the boundary complexes of cross-polytopes $C_d^{\Delta}$,
\item the boundary complexes of the McMullen-Shephard polytopes $H_{\,2d+2}^{\,d}$
      for even $d$,
\item and the spheres $C\!S^{\,3}_{\,4m}$ for $m\geq 3$.
\end{itemize}
In addition, the multiple join product $(S^d)^{*r}$ and the wreath
product $\partial\Delta_r\wr S^d$ of Joswig and Lutz \cite{JoswigLutz2003pre} 
provide two constructions to obtain \emph{derived series} of
vertex-transitive spheres for every vertex-transitive simplicial sphere $S^d$.
This way, it is even possible to get series of vertex-transitive
non-PL spheres \cite{JoswigLutz2003pre}.

The boundaries of tricyclic or multicyclic polytopes might yield further
series of vertex-transitive spheres, but it is seemingly a difficult problem 
to determine for which parameters these polytopes are simplicial. 
(Three examples of simplicial tricyclic $6$-polytopes were identified in \cite[Ch.~2]{Lutz1999}.)

Various series of vertex-transitive triangulations of surfaces
can be found in the literature; see, for example,
\cite{Altshuler1971}, \cite{Heffter1891},
\cite{KuehnelLassmann1996-bundle}, and \cite{Ringel1961}.

In higher dimensions, however, we know, apart from
the above vertex-tran\-si\-tive spheres, of only one additional
three-parameter family $M^d_k(n)$ of vertex-tran\-si\-tive triangulations 
due to K\"uhnel and Lassmann \cite{KuehnelLassmann1996-bundle}.
The combinatorial manifolds $M^d_k(n)$ on $n\geq 2^{d-k}(k+3)-1$ vertices for
$k=1,\dots,d-1$ are $k$-sphere bundles over the $(d-k)$-dimensional torus
and are invariant under the standard vertex-transitive action of the dihedral group $D_n$.
In particular, $M^d_1(n)$ is a vertex-transitive triangulation of the $d$-dimensional
torus with $n\geq 2^{d+1}-1$ vertices, and, as an additional case,
$M^d_d(d+2)$ is the boundary of the $(d+1)$-simplex; see also \cite{Kuehnel1986a-series},
\cite{KuehnelLassmann1985-di}, and \cite{KuehnelLassmann1988-dtori}.

\subsection*{4\,\, Products of Spheres}
\label{sec:products_of_spheres}

The following inequalities hold for centrally symmetric
combinatorial $2$- and $4$-manifolds $M$
with Euler characteristic $\chi(M)$.

\begin{thm} {\rm (K\"uhnel \cite{Kuehnel1996-central})}
Let $M$ be a centrally symmetric surface with\linebreak
$n=2k$\, vertices. 
Then
  \begin{equation}
    -3\,(\chi (M)-2)\leq 4^2\binom{\frac{1}{2}(k-1)}{2},
  \end{equation}
with equality if and only if $M$ contains the $1$-skeleton
of the $k$-dimensional crosspolytope $C_k^{\Delta}$, i.e.,
if $M$ is centrally $2$-neighborly.
\end{thm}

\begin{thm} {\rm (Sparla \cite[4.8]{Sparla1997}, \cite{Sparla1998})}
Let $M$ be a centrally symmetric combinatorial $4$-manifold with\, $n=2k$\, vertices.
Then
  \begin{equation}
    10\,(\chi (M)-2)\leq 4^3\binom{\frac{1}{2}(k-1)}{3},
  \end{equation}
with equality if and only if $M$ contains the $2$-skeleton
of the $k$-dimensional crosspolytope $\partial C_k^{\Delta}$, i.e.,
if $M$ is centrally $3$-neighborly.
\end{thm}

There are essentially two ways to make use of these bounds.
For fixed number $n=2k$ of vertices they give restrictions
on the Euler characteristic $\chi(M)$ of a centrally
symmetric combinatorial $2$- respectively $4$-manifold $M$ with $n$ vertices.
On the other hand, they provide \emph{lower bounds} 
on the number of vertices $n$ of a centrally
symmetric combinatorial $2$- respectively $4$-manifold $M$ 
with \emph{given} Euler characteristic $\chi(M)$.

Sparla conjectured a generalization of these bounds 
to centrally symmetric combinatorial $2r$-manifolds.

\begin{conj} {\rm (Sparla \cite[4.11]{Sparla1997}, \cite{Sparla1998})}
Let $M$ be a centrally symmetric combinatorial $2r$-man\-i\-fold with\, $n=2k$\, vertices.
Then
  \begin{equation}\label{eq:sparla_conj}
    (-1)^{r}\binom{2r+1}{r+1}(\chi (M)-2) \leq 4^{r+1}\binom{\frac{1}{2}(k-1)}{r+1},
  \end{equation}
with equality if and only if $M$ contains the $r$-skeleton
of the $k$-dimensional crosspolytope $\partial C_k^{\Delta}$, i.e.,
if $M$ is centrally $(r+1)$-neighborly.
\end{conj}

Sparla's conjecture is known to hold for $r=1$ and $r=2$ (see above) as well as in the following cases (cf.\ \cite{Novik2003pre} and \cite[4.12]{Sparla1997}):
\begin{itemize}
\item $n=4r+2$, where we trivially have $M =\partial C_{2r+1}^{\Delta}$,
\item $n\geq 4r+4$\,\, and\, $\left\{ \begin{array}{ll}\chi(M)\leq 2 & \quad\mbox{if $r$ is even,}\\
                                                       \chi(M)\geq 2 & \quad\mbox{if $r$ is odd,} \end{array}  \right.$
\item $n\geq 6r+3$\,\, (Novik \cite{Novik2003pre}).
\end{itemize}

For the sphere products $S^r\!\times\!S^r$
we have $(-1)^{r}(\chi (S^r\!\times\!S^r)-2)=2$, 
since $\chi(S^r\!\times\!S^r)=4$ if $r$ is even
and   $\chi(S^r\!\times\!S^r)=0$ if $r$ is odd.
In particular, for $n=4r+4$, i.e., for $k=2r+2$, the inequality (\ref{eq:sparla_conj}) becomes
equality, $2\binom{2r+1}{r+1}= 4^{r+1}\binom{\frac{1}{2}(2r+1)}{r+1}$ (see \cite[p.~70]{Sparla1997}).
Therefore, Sparla's conjecture, if true, would imply that
centrally symmetric combinatorial triangulations of the sphere products 
$S^r\!\times\!S^r$ with $4r+4$ vertices must contain the $r$-skeleton of $\partial C_{2r+2}^{\Delta}$.

\begin{conj} {\rm (Sparla \cite{Sparla1998})}
There are centrally $(r+1)$-neighborly trian\-gulations
of the sphere products $S^r\!\times\!S^r$ on $4r+4$ vertices.
\end{conj}

A centrally $2$-neighborly triangulation of the $2$-torus
with $8$ vertices is well known (cf.\ \cite[$\mbox{}^2\hspace{.3pt}8^{\,15}_{\,1}$]{Lutz1999}).
Centrally $3$-neighborly triangulations of the product $S^2\!\times\!S^2$ were first found by 
Sparla~\cite{Sparla1997} and by Lassmann and Sparla~\cite{LassmannSparla2000}:
There are precisely three centrally $3$-neighborly triangulations of 
$S^2\!\times\!S^2$ with $12$ vertices that have a vertex-transitive cyclic group action.


Our search for nearly neighborly centrally symmetric spheres  
with the program MANIFOLD\_VT also produced centrally symmetric 
triangulations of $d$-dimensional products of sphe\-res with $n=2d+4$ vertices,
denoted by the symbols $\mbox{}^d_{\times}\hspace{.3pt}n^{\,di/cy}_{\,z}$.
In fact, we completely enumerated
all such manifolds with a vertex-transitive cyclic or dihedral group
action for the parameters listed in Table~\ref{tbl:cent-prod}.
For 8-manifolds with 20 vertices, an enumeration was only possible
for the dihedral group action.

\begin{thm}\label{thm:cent-products}
For the products of spheres
$$\begin{array}{lllllll}
S^1\!\times\!S^1, & S^2\!\times\!S^1, & S^3\!\times\!S^1, & S^4\!\times\!S^1, & S^5\!\times\!S^1, & S^6\!\times\!S^1, & S^7 \!\times\!S^1, \\
                  &                   & S^2\!\times\!S^2, & S^3\!\times\!S^2, &                   & S^5\!\times\!S^2, &  \\
                  &                   &                   &                   & S^3\!\times\!S^3, & S^4\!\times\!S^3, & S^5\!\times\!S^3, \\
                  &                   &                   &                   &                   &                   & S^4\!\times\!S^4
\end{array}$$
there are centrally symmetric (combinatorial) triangulations
with a vertex-transitive dihedral group action on\, $n = 2d+4$\, vertices.
However, there is no sphere product\, $S^4\!\times\!S^2$\,
with a vertex-transitive cyclic group action on $16$ vertices
and no sphere product\, $S^6\!\times\!S^2$\,
with a vertex-transitive dihedral group action on $20$ vertices.
\end{thm}

\noindent
\textbf{Proof.}
The examples of Theorem~\ref{thm:cent-products} are listed in Table~\ref{tbl:cent-prod}.
We used the program BISTELLAR \cite{Lutz_BISTELLAR} to verify that 
in each case the link of vertex $1$ and therefore, by
vertex-transitivity, all vertex-links are combinatorial spheres. 
Hence, the examples are combinatorial manifolds.
The homology of the manifolds was computed 
with the program HOMOLOGY by Heckenbach~\cite{Heckenbach1997}
and, in each case, is that of a product of spheres.

The topological types of the examples $S^{d-1}\!\times\!S^1$ were determined in \cite{KuehnelLassmann1996-bundle},
and Sparla~\cite{Sparla1997} showed that the examples $\mbox{}^4_{\times}\hspace{.3pt}12^{\,cy}_{\,1}$,
$\mbox{}^4_{\times}\hspace{.3pt}12^{\,cy}_{\,2}$, and $\mbox{}^4_{\times}\hspace{.3pt}12^{\,di}_{\,1}$
are triangulations of $S^2\!\times\!S^2$. 
All remaining examples are simply connected, since they are at least
centrally $3$-neighborly. Each $d$-dimensional example occurs as a
subcomplex of the $(d+1)$-dimensional boundary sphere $\partial C_{d+2}^{\Delta}$
of the crosspolytope $C_{d+2}^{\Delta}$. According to Kreck \cite{Kreck2001}
every simply connected $d$-dimensional submanifold of the sphere $S^{d+1}$ 
with the homology of $S^{d-r}\!\times\!S^r$, $1<r\leq d/2$, 
is homeomorphic to\, $S^{d-r}\!\times\!S^r$. 
Therefore, all the examples of Table~\ref{tbl:cent-prod} 
are products of spheres.
\hfill $\Box$

\begin{conj}
There is a centrally $(\lfloor\frac{d}{2}\rfloor +1)$-neighborly
(combinatorial) \linebreak
triangulation of
every product of spheres\, $S^{\lceil\frac{d}{2}\rceil}\!\times\!S^{\lfloor\frac{d}{2}\rfloor}$\,
with a vertex-transitive dihedral group action on\, $n = 2d+4$\, vertices.
\end{conj}

\begin{landscape}

\small
\defaultaddspace=.1em

\setlength{\LTleft}{0pt}
\setlength{\LTright}{0pt}
\begin{longtable}{@{\extracolsep{\fill}}rrlllll@{}}
\caption{\protect\parbox[t]{15cm}{Centrally symmetric products of spheres with
  $n\!=\!2d\!+\!4$ vertices and cyclic group action.}}\label{tbl:cent-prod}
\\\toprule
 $d$ & $n$ & Manifold &  $f$-vector   & Type                             & List of orbits & Remarks \\ \midrule
\endfirsthead
\caption{\protect\parbox[t]{15cm}{Centrally symmetric products of spheres (continued).}}
\\\toprule
 $d$ & $n$ &   $f$-vector   & Type                              & List of orbits & Remarks\\ \midrule
\endhead
\bottomrule
\endfoot
 \addlinespace
 \addlinespace
 \addlinespace
 \addlinespace
  2  &  8  &      $S^1\!\times\!S^1$ &            (\emph{24},16) & $\mbox{}^2_{\times}\hspace{.3pt}8^{\,di}_{\,1}$ & $123_{\,8}$\, $136_{\,8}$ & \cite[$M^2_1(8)$]{KuehnelLassmann1996-bundle}, \\
     &     &                         &                          &        && \cite[$\mbox{}^2\hspace{.3pt}8^{\,15}_{\,1}$]{Lutz1999} \\
 \addlinespace
  3  & 10  &      $S^2\!\times\!S^1$ &               (\emph{40},60,30) & $\mbox{}^3_{\times}\hspace{.3pt}10^{\,di}_{\,1}$& $1235_{\,20}$\, $1245_{\,10}$ & \cite{Walkup1970}, \\
     &     &                         &                          &        && \cite[$M^3_2(10)$]{KuehnelLassmann1996-bundle}, \\
     &     &                         &                          &        && \cite[$\mbox{}^3\hspace{.3pt}10^{\,3}_{\,2}$]{Lutz1999} \\
 \addlinespace
  4  & 12  &      $S^3\!\times\!S^1$ &          (\emph{60},120,120,48) & $\mbox{}^4_{\times}\hspace{.3pt}12^{\,di}_{\,2}$& $12346_{\,24}$\, $12356_{\,24}$ & \cite[$M^4_3(12)$]{KuehnelLassmann1996-bundle}, \\
     &     &                         &                          &        && \cite[$\mbox{}^4\hspace{.3pt}12^{\,12}_{\,1}$]{Lutz1999} \\
 \addlinespace
     &     &      $S^2\!\times\!S^2$ &          (\emph{60},\emph{160},180,72) & $\mbox{}^4_{\times}\hspace{.3pt}12^{\,cy}_{\,1}$& $12345_{\,12}$\, $12356_{\,12}$\, $1236\,11_{\,12}$\, $12569_{\,12}$\, $1269\,11_{\,12}$\, $1358\,10_{\,12}$ & \cite[$M_1$]{Sparla1997}, \\
     &     &                         &                          &        && \cite[$\mbox{}^4\hspace{.3pt}12^{\,11}_{\,1}$]{Lutz1999} \\
 \addlinespace
     &     &                         &                          & $\mbox{}^4_{\times}\hspace{.3pt}12^{\,cy}_{\,2}$& $12345_{\,12}$\, $12356_{\,12}$\, $1236\,11_{\,12}$\, $1256\,10_{\,12}$\, $1269\,11_{\,12}$\, $1358\,10_{\,12}$ & \cite[$M\!=\!M_2$]{Sparla1997}, \\
     &     &                         &                          &        && \cite{Sparla1998}, \\
     &     &                         &                          &        && \cite[$\mbox{}^4\hspace{.3pt}12^{\,124}_{\,1}$]{Lutz1999} \\
 \addlinespace
     &     &                         &                          & $\mbox{}^4_{\times}\hspace{.3pt}12^{\,di}_{\,1}$& $12345_{\,12}$\, $1235\,10_{\,24}$\, $1236\,10_{\,12}$\, $12459_{\,12}$\, $1358\,10_{\,12}$ & \cite[$M_3$]{Sparla1997}, \\
     &     &                         &                          &        && \cite[$\mbox{}^4\hspace{.3pt}12^{\,28}_{\,1}$]{Lutz1999} \\
 \addlinespace
  5  & 14  &      $S^4\!\times\!S^1$ &      (\emph{84},210,280,  & $\mbox{}^5_{\times}\hspace{.3pt}14^{\,di}_{\,1}$& $123457_{\,28}$\, $123467_{\,28}$\, $123567_{\,14}$ & \cite[$M^5_4(14)$]{KuehnelLassmann1996-bundle} \\
     &     &                         &      210,70)             &        && \cite[$\mbox{}^5\hspace{.3pt}14^{\,3}_{\,8}$]{Lutz1999} \\
 \addlinespace
     &     &      $S^3\!\times\!S^2$ &     (\emph{84},\emph{280},490, & $\mbox{}^5_{\times}\hspace{.3pt}14^{\,di}_{\,2}$& $123467_{\,28}$\, $12346\,12_{\,28}$\, $123567_{\,14}$\, $12357\,11_{\,28}$ & \cite[$\mbox{}^5\hspace{.3pt}14^{\,3}_{\,9}$]{Lutz1999} \\ 
     &     &                         &      420,140)            &        & $12457\,13_{\,14}$\, $1246\,10\,12_{\,28}$ & \\
 \addlinespace
  6  & 16  &      $S^5\!\times\!S^1$ & (\emph{112},336,560, & $\mbox{}^6_{\times}\hspace{.3pt}16^{\,di}_{\,2}$& $1234568_{\,32}$\, $1234578_{\,32}$\, $1234678_{\,32}$ & \cite[$M^6_5(16)$]{KuehnelLassmann1996-bundle} \\
     &     &                         & 560,336,96) &&&       \\
 \addlinespace
\newpage
 \addlinespace
 \addlinespace
 \addlinespace
  \addlinespace
     &     &      $S^3\!\times\!S^3$ & (\emph{112},\emph{448},\emph{1120},&$\mbox{}^6_{\times}\hspace{.3pt}16^{\,cy}_{\,1}$& $1234567_{\,16}$\, $1234578_{\,16}$\, $123458\,15_{\,16}$\, $123478\,13_{\,16}$ & \\ 
     &     &                         & 1568,1120,320)           && $12347\,13\,14_{\,16}$\, $12348\,13\,15_{\,16}$\, $1234\,13\,14\,15_{\,16}$\, $123568\,15_{\,16}$ & \\ 
     &     &                         &                          && $123678\,12_{\,16}$\, $12368\,12\,13_{\,16}$\, $12368\,13\,15_{\,16}$\, $12378\,12\,13_{\,16}$ & \\ 
     &     &                         &                          && $124578\,11_{\,16}$\, $12457\,11\,14_{\,16}$\, $12458\,11\,14_{\,16}$\, $12478\,11\,13_{\,16}$ & \\
     &     &                         &                          && $1247\,11\,13\,14_{\,16}$\, $1248\,11\,13\,15_{\,16}$\, $1268\,11\,13\,15_{\,16}$\, $1357\,10\,12\,14_{\,16}$ & \\
 \addlinespace
     &     &                         &                          & $\mbox{}^6_{\times}\hspace{.3pt}16^{\,di}_{\,1}$& $1234567_{\,16}$\, $1234578_{\,32}$\, $123458\,14_{\,16}$\, $123478\,13_{\,32}$ & \\ 
     &     &                         &                          && $123567\,12_{\,16}$\, $12356\,12\,15_{\,32}$\, $123578\,12_{\,32}$\, $12358\,12\,15_{\,16}$ & \\ 
     &     &                         &                          && $12378\,12\,13_{\,16}$\, $12458\,11\,14_{\,16}$\, $12467\,11\,13_{\,16}$\, $12468\,11\,13_{\,32}$ & \\
     &     &                         &                          && $1247\,11\,13\,14_{\,32}$\, $1357\,10\,12\,14_{\,16}$ & \\
 \addlinespace
  7  & 18  &      $S^6\!\times\!S^1$ & (\emph{144},504,1008, & $\mbox{}^7_{\times}\hspace{.3pt}18^{\,di}_{\,1}$& $12345679_{\,36}$\, $12345689_{\,36}$\, $12345789_{\,36}$\, $12346789_{\,18}$ & \cite[$M^7_6(18)$]{KuehnelLassmann1996-bundle} \\
     &     &                         & 1260,1008, &&&       \\
     &     &                         & 504,126) &&&       \\
 \addlinespace
     &     &      $S^5\!\times\!S^2$ & (\emph{144},\emph{672},1764, & $\mbox{}^7_{\times}\hspace{.3pt}18^{\,di}_{\,2}$& $12345689_{\,36}$\, $1234568\,16_{\,36}$\, $12345789_{\,36}$\, $1234579\,15_{\,36}$ & \\
     &     &                         & 2772,2688,                 && $12346789_{\,18}$\, $1234679\,17_{\,36}$\, $123468\,14\,16_{\,36}$\, $1235679\,17_{\,18}$ & \\
     &     &                         & 1512,378)                  && $123579\,13\,15_{\,36}$\, $123579\,13\,17_{\,36}$\, $124579\,15\,17_{\,18}$\, $12468\,12\,14\,16_{\,36}$ & \\
 \addlinespace
     &     &      $S^4\!\times\!S^3$ & (\emph{144},\emph{672},\emph{2016}, & $\mbox{}^7_{\times}\hspace{.3pt}18^{\,di}_{\,3}$& $1234579\,15_{\,36}$\, $1234579\,17_{\,36}$\, $1234679\,14_{\,36}$\, $1234679\,17_{\,36}$ & \\
     &     &                         & 3780,4200,               && $123467\,14\,17_{\,36}$\, $1234689\,14_{\,36}$\, $1234689\,16_{\,36}$\, $123479\,14\,15_{\,36}$ & \\
     &     &                         & 2520,630)                && $1235679\,13_{\,36}$\, $1235679\,17_{\,18}$\, $1235689\,13_{\,36}$\, $1235689\,16_{\,36}$ & \\
     &     &                         &                          && $123569\,16\,17_{\,36}$\, $1235789\,13_{\,36}$\, $123589\,15\,16_{\,36}$\, $123679\,13\,14_{\,36}$ & \\
     &     &                         &                          && $123689\,13\,14_{\,36}$\, $123789\,13\,15_{\,18}$\, $124589\,15\,16_{\,18}$ & \\
 \addlinespace
  8  & 20  &      $S^7\!\times\!S^1$ & (\emph{180},720,1680, & $\mbox{}^8_{\times}\hspace{.3pt}20^{\,di}_{\,2}$& $12345678\,10_{\,40}$\, $12345679\,10_{\,40}$\, $12345689\,10_{\,40}$\, $12345789\,10_{\,40}$ & \cite[$M^8_7(20)$]{KuehnelLassmann1996-bundle} \\
     &     &                         & 2520,2520,1680, &&&       \\
     &     &                         & 720,160) &&&       \\
 \addlinespace
\newpage
 \addlinespace
 \addlinespace
 \addlinespace
  \addlinespace
     &     &      $S^5\!\times\!S^3$ & (\emph{180},\emph{960},\emph{3360}, & $\mbox{}^8_{\times}\hspace{.3pt}20^{\,di}_{\,3}$& $12345689\,10_{\,40}$\, $12345689\,17_{\,40}$\, $1234568\,10\,19_{\,40}$\, $1234569\,10\,18_{\,40}$ & \\
     &     &                         & 7560,10920,9840,         && $12345789\,10_{\,40}$\, $1234578\,10\,16_{\,40}$\, $1234579\,10\,18_{\,40}$ & \\
     &     &                         & 5040,1120)               && $123457\,10\,16\,18_{\,40}$\, $1234678\,10\,19_{\,40}$\, $1234679\,10\,15_{\,40}$ & \\
     &     &                         &                          && $123467\,10\,15\,19_{\,40}$\, $123469\,10\,15\,18_{\,40}$\, $123469\,15\,17\,18_{\,40}$ & \\
     &     &                         &                          && $123479\,10\,15\,18_{\,40}$\, $1235679\,10\,18_{\,40}$\, $1235689\,14\,17_{\,40}$ & \\
     &     &                         &                          && $123568\,14\,17\,19_{\,40}$\, $123569\,14\,17\,18_{\,40}$\, $123578\,10\,16\,19_{\,40}$ & \\
     &     &                         &                          && $123679\,10\,15\,18_{\,40}$\, $124579\,10\,16\,18_{\,40}$\, $124579\,13\,16\,18_{\,40}$ & \\
     &     &                         &                          && $12458\,10\,13\,16\,17_{\,40}$\, $124679\,13\,15\,18_{\,40}$\, $12469\,10\,13\,15\,18_{\,40}$ & \\
     &     &                         &                          && $12478\,10\,13\,15\,16_{\,40}$\, $12478\,10\,13\,15\,19_{\,40}$\, $12478\,10\,15\,16\,19_{\,40}$ & \\
 \addlinespace
     &     &      $S^4\!\times\!S^4$ & (\emph{180},\emph{960},\emph{3360}, & $\mbox{}^8_{\times}\hspace{.3pt}20^{\,di}_{\,1}$& $123456789_{\,20}$\, $12345679\,10_{\,40}$\, $1234567\,10\,18_{\,20}$\, $1234569\,10\,17_{\,40}$ & \\
     &     &                         & \emph{8064},12600,12000,  && $12345789\,16_{\,40}$\, $1234578\,16\,19_{\,40}$\, $1234579\,10\,16_{\,40}$ & \\
     &     &                         & 6300,1400)               && $123457\,10\,16\,19_{\,20}$\, $123459\,10\,16\,17_{\,20}$\, $12346789\,15_{\,20}$ & \\
     &     &                         &                          && $1234679\,10\,18_{\,40}$\, $1234679\,15\,18_{\,40}$\, $123469\,10\,17\,18_{\,40}$ & \\
     &     &                         &                          && $1234789\,15\,16_{\,40}$\, $123479\,10\,16\,18_{\,40}$\, $1235679\,10\,14_{\,40}$ & \\

     &     &                         &                          && $123567\,10\,14\,18_{\,20}$\, $1235689\,10\,14_{\,40}$\, $123568\,10\,14\,17_{\,40}$ & \\
     &     &                         &                          && $12356\,10\,14\,17\,19_{\,40}$\, $12356\,10\,14\,18\,19_{\,20}$\, $1235789\,14\,16_{\,20}$ & \\
     &     &                         &                          && $123578\,14\,16\,19_{\,40}$\, $123579\,10\,14\,16_{\,40}$\, $12357\,10\,14\,16\,19_{\,40}$ & \\
     &     &                         &                          && $123589\,10\,14\,16_{\,40}$\, $123589\,10\,16\,17_{\,40}$\, $12358\,10\,14\,16\,19_{\,20}$ & \\
     &     &                         &                          && $123679\,10\,14\,18_{\,40}$\, $123679\,14\,15\,18_{\,40}$\, $12368\,10\,14\,15\,17_{\,40}$ & \\
     &     &                         &                          && $12369\,10\,14\,15\,18_{\,20}$, $123789\,14\,15\,16_{\,20}$\, $124579\,10\,13\,16_{\,40}$ & \\
     &     &                         &                          && $12457\,10\,13\,16\,19_{\,20}$\, $12459\,10\,13\,16\,17_{\,20}$\, $12467\,10\,13\,15\,18_{\,20}$ & \\
     &     &                         &                          && $124689\,13\,15\,17_{\,20}$\, $12468\,10\,13\,15\,17_{\,40}$\, $12469\,10\,15\,17\,18_{\,40}$ & \\
     &     &                         &                          && $12469\,13\,15\,17\,18_{\,40}$\, $12479\,10\,13\,15\,18_{\,40}$\, $13579\,12\,14\,16\,18_{\,20}$ & \\
 \addlinespace
     &     &         ........        &         ........         & $\mbox{}^8_{\times}\hspace{.3pt}20^{\,cy}_{\,...}$ & ........ & \\
 \addlinespace
\end{longtable}
\end{landscape}

\bibliography{}

\bigskip
\medskip

\noindent
Frank H.\ Lutz\\
Technische Universit\"at Berlin\\
Fakult\"at II - Mathematik und Naturwissenschaften\\
Institut f\"ur Mathematik, Sekr. MA 6-2\\
Stra\ss e des 17.\ Juni 136\\
D-10623 Berlin\\
{\tt lutz@math.tu-berlin.de}

\end{document}